\documentclass[10pt]{article}
\usepackage{amssymb,amsmath,amsthm, amsfonts}

\def\sqr#1#2{{\vcenter{\vbox{\hrule height.#2pt
              \hbox{\vrule width.#2pt height#1pt \kern#1pt \vrule width.#2pt}
          \hrule height.#2pt}}}}
\def\signed #1{{\unskip\nobreak\hfil\penalty50
          \hskip2em\hbox{}\nobreak\hfil#1
          \parfillskip=0pt \finalhyphendemerits=0 \par}}
\def\endpf{\signed {$\sqr69$}}
\def\sqr#1#2{{\vcenter{\vbox{\hrule height.#2pt
              \hbox{\vrule width.#2pt height#1pt \kern#1pt \vrule width.#2pt}
              \hrule height.#2pt}}}}
\def\signed #1{{\unskip\nobreak\hfil\penalty50
              \hskip2em\hbox{}\nobreak\hfil#1
              \parfillskip=0pt \finalhyphendemerits=0 \par}}
\def\endpf{\signed {$\sqr69$}}
\def\3n{\negthinspace \negthinspace \negthinspace }
\def\2n{\negthinspace \negthinspace }
\def\1n{\negthinspace }

\def\={\buildrel \triangle \over =}

\def\O{\Omega}

\def\q{\quad}
\def\qq{\qquad}

\def\max{\mathop{\rm max}}

\def\exp{\mathop{\rm exp}}
\def\sup{\mathop{\rm sup}}

\def\|{\Big |}
\def\({\Big (}
\def\){\Big )}
\def\[{\Big[}
\def\]{\Big]}
\def\be{\begin{equation}}
\def\bel{\begin{equation}\label}
\def\ee{\end{equation}}
\def\bt{\begin{theorem}}
\def\bcd{\begin{condition}}
\def\ecd{\end{condition}}
\def\et{\end{theorem}}
\def\bc{\begin{corollary}}
\def\ec{\end{corollary}}
\def\bde{\begin{definition}}
\def\ede{\end{definition}}
\def\bl{\begin{lemma}}
\def\el{\end{lemma}}
\def\bp{\begin{proposition}}
\def\ep{\end{proposition}}
\def\br{\begin{remark}}
\def\er{\end{remark}}
\def\ba{\begin{array}}
\def\ea{\end{array}}
\def\ed{\end{document}}

\def\square#1{\vbox{\hrule\hbox{\vrule height#1%
     \kern#1\vrule}\hrule}}
\def\rectangle#1#2{\vbox{\hrule\hbox{\vrule height#1%
     \kern#2\vrule}\hrule}}

\font\tenbb=msbm10 \font\sevenbb=msbm7 \font\fivebb=msbm5

\newfam\bbfam
\scriptscriptfont\bbfam=\fivebb \textfont\bbfam=\tenbb
\scriptfont\bbfam=\sevenbb

\newtheorem{lemma}{Lemma}[section]
\newtheorem{remark}{Remark}[section]

\newtheorem{theorem}{Theorem}[section]
\newtheorem{corollary}{Corollary}[section]

\newtheorem{definition}{Definition}[section]
\newtheorem{proposition}{Proposition}[section]
\newtheorem{condition}{Condition}[section]

\makeatletter
   
   \@addtoreset{equation}{section}
\makeatother

\begin{document}

\title{Integral-Partial Differential Equations
of Isaacs' Type Related to Stochastic Differential Games with
Jumps}

\author{Rainer Buckdahn{\footnote{This work has been done in the frame of the Marie Curie ITN Project ``Deterministic and Stochastic Controlled Systems and Applications", call: F97-PEOPLE-2007-1-1-ITN, no: 213841-2. \qq \qq \qq \qq \qq  \mbox{ }  \mbox{ } \mbox{ }  \mbox{ } \mbox{ }  \mbox{ }\mbox{ }  \mbox{ }\mbox{ }  \mbox{ } \mbox{   } \mbox{   } \mbox{ } \qq \qq \qq \qq \qq \qq \qq \qq \qq \qq
$ ^{\ast\ast}$Corresponding author. The work of Juan LI has been supported by a
one-year fellowship awarded by the General Council of Finist\`{e}re,
France, and by the NSF of P.R.China (No. 10701050), Shandong
Province (No. Q2007A04), Independent Innovation Foundation of Shandong University, SRF for ROCS (SEM) and National Basic Research Program of China (973 Program) (No. 2007CB814904).}}\\
{\small D\'{e}partement de Math\'{e}matiques, Universit\'{e} de
Bretagne Occidentale,}\\
 {\small 6, avenue Victor-le-Gorgeu, CS 93837, 29238 Brest
Cedex 3, France.}\\
{\small{\it E-mail: Rainer.Buckdahn@univ-brest.fr.}}\\
Ying Hu$ ^{\ast}$\\
{\small IRMAR, Universit\'{e} Rennes 1, Campus de Beaulieu, 35042 Rennes Cedex, France.}\\
{\small{\it ying.hu@univ-rennes1.fr.}}\\
 Juan Li$^{\ast\ast}$\\
{\small School of Mathematics and Statistics, Shandong University at Weihai, Weihai 264200, P. R. China.}\\
{\small {\it E-mail: juanli@sdu.edu.cn.}} \date{April 15, 2010}}
\maketitle

 \noindent{\bf Abstract}\hskip4mm In this paper we study
zero-sum two-player stochastic differential games with jumps with
the help of theory of Backward Stochastic Differential Equations
(BSDEs). We generalize the results of Fleming and Souganidis~\cite{FS} and those by Biswas~\cite{BIS} by considering a controlled stochastic system
driven by a d-dimensional Brownian motion and a Poisson random
measure and by associating nonlinear cost functionals defined by controlled
BSDEs. Moreover, unlike the both papers cited above we allow the admissible control processes of
both players to depend on all events occurring before the beginning
of the game. This quite natural extension allows the players
to take into account such earlier events, and it makes  even easier to derive the dynamic programming principle.
The price to pay is that the cost functionals become random
variables and so also the upper and the lower value functions of the game are a priori random fields. The use of a new method allows to prove that, in fact, the upper and the lower value functions
are deterministic. On the other hand, the application of BSDE
methods~\cite{P1} allows to prove a dynamic programming principle for the upper and
the lower value functions in a very straight-forward
way, as well as the fact that they are the unique
viscosity solutions of the upper and the lower integral-partial
differential equations of Hamilton-Jacobi-Bellman-Isaacs' type,
respectively. Finally, the existence of the
value of the game is got in this more general setting if Isaacs' condition
holds.
\newpage
\section{\large{Introduction}}
\hskip1cm In the present work we investigate two-player zero-sum
stochastic differential games in the framework of Brownian motion
and Poisson random measure. Fleming and Souganidis~\cite{FS} were
the first to study in a rigorous manner two-player zero-sum
stochastic differential games (SDGs). They proved that the lower and the
upper value functions of such games satisfy the dynamic programming
principle (DPP), that they are the unique viscosity solutions of the
associated Bellman-Isaacs equations and that they coincide under the Isaacs'
condition. Their work has translated former results by Evans and
Souganidis~\cite{ES} from a deterministic into the stochastic
framework and has given an important impulse for the research in the
theory of stochastic differential games. Various recent works on SDGs are based on the ideas developed in~\cite{FS}, see, for
instance, Buckdahn, Cardaliaguet, Rainer~\cite{BCR}, Hou,
Tang~\cite{HT}, Mataramvura, $\varnothing$ksental~\cite{MO} and so on, we shall, in particular, also mention the recent work by Biswas~\cite{BIS} on SDGs with jumps. The reader interested in this subject is also referred to the references
given in~\cite{FS}.

In the present work we study an extension of the results of the pioneering work of Fleming and Souganidis~\cite{FS} and those by Biswas~\cite{BIS} on SDGs with jumps. More precisely, inspired by~\cite{BL}, we consider SDGs with jumps on the Wiener-Poisson space and we allow the admissible
control processes to depend on the full past of the trajectories of
the driving Brownian motion and the Poisson random measure. This means,
in particular, that they can also depend on information occurring before
the beginning of the game. This approach combined with the notion of stochastic backward semigroups, introduced by Peng~\cite{P1}, simplifies the proof of the DPP considerably. But it also has the
consequence that the cost functionals become random variables. In~\cite{BL}, for SDGs in a Brownian setting without jumps, the authors introduced the method of Girsanov transformation, in order to prove that in spite of the randomness of the cost functionals the lower and the upper value functions of the game are deterministic. However, this method doesn't apply to SDGs with jumps. For this reason we study a new type of transformation on the Wiener-Poisson space, which allows to show that also in the case of SDGs with jumps the upper and the lower value functions are deterministic, in spite of controls which can depend on the whole past. Another extension concerns the cost functionals. We consider nonlinear ones, defined through a doubly controlled backward stochastic differential equation (BSDE) with jumps. These both extensions of the framework
in~\cite{FS},~\cite{BIS} and~\cite{BL} are crucial because they allow to
harmonize the setting for stochastic differential games  with jumps with that
for the stochastic control theory and to simplify considerably the
approach in~\cite{FS} and~\cite{BIS} by using BSDE methods.

BSDEs in the framework of Brownian motion in their general
non-linear form were introduced by Pardoux and Peng~\cite{PP1} in
1990. They have been studied since then by a lot of authors and have
found various applications, namely in stochastic control, finance
and the second order PDE theory. BSDE methods, originally developed
by Peng~\cite{P1} and ~\cite{P2} for the stochastic control theory, have
been introduced to the theory of stochastic differential games by
Hamad\`{e}ne and Lepeltier~\cite{HL} and Hamad\`{e}ne, Lepeltier and
Peng~\cite{HLP}, in order to study games with a dynamics whose diffusion
coefficient is strictly elliptic and does not depend on the controls.
BSDEs in the framework of Brownian motion and Poisson random measure
were first considered by Tang and Li~\cite{TL}, later by Barles, Buckdahn and Pardoux~\cite{BBP}, and so on. In Li and Peng~\cite{LP} they studied the stochastic control theory for BSDE with jumps.

In the present paper we study the general framework of SDGs with jumps. The dynamics of the stochastic differential game in the framework of
Brownian motion and compensated Poisson random measure we investigate is given
by a doubly controlled system of stochastic differential equations (see equation
(3.1)). The cost functionals (interpreted as a payoff for Player I
and as a cost for Player II)(see (3.7)) are introduced by a BSDE governed by a Brownian motion and a compensated Poisson random measure (see equation (3.5)). It is well known in the theory of differential games, that players cannot restrict to play only control processes: one player has to fix a strategy
while the other player chooses the best answer to this strategy in
the form of a control process. So the lower value function $W$ is
defined as the essential infimum  of the essential supremum of all
cost functionals, where the essential infimum is taken over all
admissible strategies of Player II and the essential supremum is
taken over all admissible controls of Player I. The upper value
function $U$ is defined by changing the roles of the both players:
as the essential supremum  of the essential infimum of all cost
functionals, where the essential supremum is taken over all
admissible strategies of Player I and the essential infimum is taken
over all admissible controls of Player II; for the precise
definitions see (3.9) and (3.10). The objective of our paper is to
investigate these lower and upper value functions.  The main results
of the paper state that $W$\ and $U$\ are deterministic (Proposition
3.1) continuous unique viscosity solutions of the associated Bellman-Isaacs equations
(Theorem 4.1), and they satisfy the DPP (Theorem 3.1).

We point out the fact that $W$ and $U$, introduced as combination of
essential infimum and essential supremum over a class of random
variables, are deterministic is far from being trivial. The method
developed by Peng~\cite{P1} (see also Theorem 6.1 of the present
paper) for value functions involving only control processes but not
strategies does not apply here since the strategies from ${\cal
A}_{t,T}$\ and ${\cal B}_{t,T}$ do not have, in general, any
continuity property. In~\cite{BL}, the authors used a new method, that of the Girsanov transformation, to solve this
difficulty for the stochastic differential games in the framework
of Brownian motion, but for the present situation-the SDGs driven by a Brownian motion and a compensated Poisson
random measure this method is not applicable anymore. To overcome this
difficulty we define a new type of measure-preserving and invertible transformations on the Wiener-Poisson space (see (3.11) and (3.12)). We
show in Lemma 3.1 that $W$\ and $U$\ are invariant under such transformations and in Lemma 3.2 we prove that the invariance of a random variable over the
Wiener-Poisson space with respect to these transformations implies that it is deterministic. We
emphasize that the proofs of the Lemmas 3.1 and 3.2 do not use BSDE methods. This makes this method also applicable to the other situations, such as standard stochastic control problems with jumps. The importance of the approach which considers control processes depending on events occurring before the beginning of the game, stems from that fact that, once proved that the upper and the lower value functions $W$\ and $U$\ are deterministic, Peng's notion of backward stochastic semigroups~\cite{P1} extended to the framework with jumps, allows to prove in a very straight-forward way the DPP and this without any approximation or technical notions ( $r$-strategies and $\pi$-controls) playing an essential role in~\cite{FS} and~\cite{BIS}. Moreover, our approach also allows to show directly with the help of the DPP that $W$\ and $U$\ are viscosity solutions of the associated Bellman-Isaacs equations.

Our paper is organized as follows. Section 2 recalls some elements
of the theory of BSDEs with jumps which will be needed in what follows. Section 3 introduces the
setting of stochastic differential game and its lower and upper
value functions $W$ and $U$, and it proves that both these  functions
are deterministic and satisfy the DPP. The proof of DPP is given in Section 6.2. In
Section 4 the DPP allows to derive with the help of Peng's BSDE method~\cite{P1} adapted to the framework of SDGs with jumps, that $W$ and $U$
are viscosity solutions of the associated Bellman-Isaacs equations. In Section 5 we prove the uniqueness
of viscosity solutions of the associated Bellman-Isaacs equations. Finally, after having characterized $W$ and $U$ as viscosity
solutions of associated Bellman-Isaacs equations we show that under
the Isaacs' condition $W$ and $U$ coincide (one says that the game
has a value). Finally, the Appendix recalls some complementary  results on FBSDEs
with jumps, to which we refer in our work.

\section{ {\large Preliminaries}}

 \hskip1cm Let us begin by introducing the setting for the stochastic
differential game we want to investigate. As underlying probability
space $(\Omega, {\cal{F}}, P)$\ we consider the completed product of
the Wiener space $(\Omega_1, {\cal{F}}_1, P_1)$\ and the Poisson
space $(\Omega_2, {\cal{F}}_2, P_2).$\ Here, $(\Omega_1,
{\cal{F}}_1, P_1)$\ is a Wiener space: $\Omega_1$\
is the set of continuous functions from ${\mathbb{R}}$\ to ${\mathbb{R}}^d$\
with value zero at 0, $\Omega_1= C_0({\mathbb{R}};{\mathbb{R}}^d)$\ endowed with the topology generated by the uniform convergence on compacts.
$ {\cal{F}}_1 $\ is the completed Borel
$\sigma$-algebra over $\Omega_1$, and
 $P_1$ the Wiener measure under which the d-dimensional coordinate
processes $B_s(\omega)=\omega_s,\ s\in {\mathbb{R}}_+,\ \omega\in \Omega_1,$\ and $B_{-s}(\omega)=\omega(-s),\ s\in {\mathbb{R}}_+,\ \omega\in \Omega_1,$\ are two independent d-dimensional Brownian motions. By $\{{\mathcal{F}}^{B}_s,\
s\geq 0\}$\ we denote the natural filtration generated by
$\{B_s\}_{s\geq 0}$\ and augmented by all $P_1$-null sets,
i.e.,
$${\mathcal{F}}^{B}_s=\sigma\{B_r, r\in (-\infty, s]\}\vee {\mathcal{N}}_{P_1}, s\geq 0. $$
We now introduce the Poisson space $(\Omega_2, {\cal{F}}_2, P_2).$\
For this, we let $E={\mathbb{R}}^l\setminus \{0\}$\ and endow the
space E with its Borel $\sigma$-field ${\mathcal{B}}(E)$. By a point
function p on E we understand a mapping $p:D_p\subset
{\mathbb{R}}\rightarrow E$, where the domain $D_p$\ is a countable
subset of the real line ${\mathbb{R}}$. The point function $p$
defines on ${\mathbb{R}}\times E$\ the counting measure
$\mu(p,dtde)$\ introduced by the relation
$$\mu(p, (s,t]\times \Delta)=\sharp\{r\in D_p\cap (s,t]: p(r)\in
\Delta\},\ \Delta\in {\mathcal{B}}(E),\ s, t\in {\mathbb{R}},\
s<t.$$ In the sequel we will often identify the point function $p$\
 with $\mu(p, .)$. Let now $\Omega_2$ denote the collection of all point functions $p$
on $E$ and ${\cal{F}}_2$ be the smallest $\sigma$-field on
$\Omega_2$ with respect to which all mappings
$p\rightarrow\mu(p,(s,t]\times \Delta),\, s,t \in {\mathbb{R}}, s<t,
\Delta\in{\cal B}(E)$ are measurable. On the measurable space
$(\Omega_2, {\cal{F}}_2)$ we consider the probability measure $P_2$
under which the canonical coordinate measure $\mu(p,dtde)$ becomes a
Poisson random measure with L\'evy measure $\lambda$. That means the
{\it compensator} $\hat{\mu}(dtde)=dt\lambda(de)$ of $\mu$
  transforms $\{\tilde{\mu}((s,t]\times
  A)=(\mu-\hat{\mu})((s,t]\times A)\}_{s\leq t}$\ to a martingale for
  any $A\in{\cal{B}}(E)$ satisfying $\lambda(A)<\infty$. Here
$\lambda$\ is an arbitrarily given $\sigma$-finite L\'evy measure on
$(E,{\cal B}(E))$, i.e., a measure on $(E,{\mathcal{B}}(E))$ with
the property that $\int_E(1\wedge|e|^2)\lambda(de)<\infty$. We
complete the probability space $(\Omega_2, {\cal{F}}_2, P_2)$ and
introduce the filtration $({\cal F}^\mu_t)_{t\geq 0}$ generated by
our coordinate measure $\mu$ by setting

$$\dot{\cal F}_t^{\mu}=\sigma\{\mu((s,r]\times \Delta):\,
-\infty<s\leq r\leq t, \Delta\in{\cal B}(E)\},\ t\geq 0,$$

\noindent and taking the right-limits ${\cal
F}_t^{\mu}=\left(\cap_{\{s>t\}}\dot{\cal F}^{\mu}_s\right)\vee{\cal
N}_{P_2},\, t\geq 0$, augmented by the $P_2$-null sets.
Finally, we put $\Omega=\Omega_1\times\Omega_2,\ {\cal F}={\cal
F}_1\otimes{\cal F}_2, P=P_1\otimes P_2,$ where ${\cal F}$ is
completed with respect to P, and the filtration ${\mathbb{F}}=\{{\cal F}_t\}_{t\geq 0}$\ is generated by
 $${\cal F}_t:={\cal F}_t^{B,\mu}={\cal F}_t^B\otimes{\cal F}_t^\mu, \ \ t\geq 0,\ \ \mbox{augmented by all P-null sets}.$$

  Let $T>0$\ be an arbitrarily fixed time horizon. For any
    $n\geq 1,$\ $|z|$ denotes the Euclidean norm of $z\in
    {\mathbb{R}}^{n}$. We  introduce also the following three spaces of processes which will be used frequently in the sequel:
\vskip0.2cm
    ${\cal{S}}^2(0, T; {\mathbb{R}}):=\{(\psi_t)_{0\leq t\leq T}\mbox{ real-valued}\ {\mathbb{F}}\mbox{-adapted c\`{a}dl\`{a}g
    process}:\\ \mbox{ }\hskip6cm
    E[\sup\limits_{0\leq t\leq T}| \psi_{t} |^2]< +\infty \}; $
    \vskip0.2cm

   ${\cal{H}}^{2}(0,T;{\mathbb{R}}^{n}):=\{(\psi_t)_{0\leq t\leq T}\ {\mathbb{R}}^{n}\mbox{-valued}\ {\mathbb{F}}
    \mbox{-progressively measurable process}:\\ \mbox{ }\hskip6cm
     \parallel\psi\parallel^2=E[\int^T_0| \psi_t| ^2dt]<+\infty \}; $
\vskip0.2cm

${\cal{K}}_{\lambda}^{2}(0,T;{\mathbb{R}}^{n}):=\{\mbox{mapping}\ K:
\Omega\times [0, T]\times E\rightarrow {\mathbb{R}}^{n}\
{\cal{P}}\otimes {\cal{B}}(E)\mbox{-measurable}:\\
\mbox{ }\hskip6cm
     \parallel K\parallel^2=E[\int^T_0\int_E| K_t(e)| ^2\lambda(de)dt]<+\infty \}. $
\footnote{${\cal{P}}$\ denotes the $\sigma$-algebra of
${\cal{F}}_t$-predictable subsets of $\Omega\times[0, T]$.}
 \vskip0.2cm
 Let us now consider a function $g:
\Omega\times[0,T]\times {\mathbb{R}} \times {\mathbb{R}}^{d}\times
L^2(E, {\cal{B}}(E), \lambda; {\mathbb{R}}) \rightarrow {\mathbb{R}}
$ with the property that $(g(t, y, z, k))_{t\in [0, T]}$ is
${\cal{P}}$-measurable for each $(y, z, k)$ in ${\mathbb{R}} \times
{\mathbb{R}}^{d}\times L^2(E, {\cal{B}}(E), \lambda; {\mathbb{R}})$,
and we also make the following assumptions on $g $ throughout the
paper:
 \vskip0.2cm

(A1) There exists a constant $C\ge 0$  such that, P-a.s., for all
$t\in [0, T],\ y_{1}, y_{2}\in {\mathbb{R}},\ z_{1}, z_{2}\in
{\mathbb{R}}^d, \ k_{1}, k_{2}\in L^2(E, {\cal{B}}(E), \lambda;
{\mathbb{R}}),\\ \mbox{ }\hskip3cm |g(t, y_{1}, z_{1}, k_{1}) - g(t,
y_{2}, z_{2}, k_{2})|\leq C(|y_{1}-y_{2}| +
|z_{1}-z_{2}|+||k_{1}-k_{2}||).$
 \vskip0.2cm

(A2) $g(\cdot,0,0,0)\in {\cal{H}}^{2}(0,T;{\mathbb{R}})$.
\vskip0.2cm

 The following result on BSDEs with jumps is by now well known, for its proof
 the reader is referred to Lemma 2.4 in Tang and Li~\cite{TL} or Theorem 2.1 in Buckdahn, Barles and Pardoux~\cite{BBP}.
 \bl Under the assumptions (A1) and (A2), for any random variable $\xi\in L^2(\O, {\cal{F}}_T,$ $P),$ the
BSDE with jump
 \be y_t = \xi + \int_t^Tg(s,y_s,z_s, k_s)ds - \int^T_tz_s\,
dB_s-\int^T_t\int_E k_s (e)\tilde{\mu} (ds,de),\q 0\le t\le T,
\label{BSDE} \ee
 has a unique adapted solution
$$(y^{T, g, \xi}_t, z^{T, g,
\xi}_t, k^{T, g, \xi}_t)_{t\in [0, T]}\in {\cal{S}}^2(0, T;
{\mathbb{R}})\times
{\cal{H}}^{2}(0,T;{\mathbb{R}}^{d})\times{\cal{K}}_{\lambda}^{2}(0,T;{\mathbb{R}}).
$$ \el
 In the sequel, we
always assume that the driving coefficient $g$\ of a BSDE with jump
satisfies (A1) and (A2).

   We  recall also the following both basic results on BSDEs with jumps.
   We begin with the well-known comparison theorem (see
   Barles, Buckdahn and Pardoux~\cite{BBP}, Proposition 2.6).

\bl (Comparison Theorem) Let $h:\Omega\times[0, T]\times
{\mathbb{R}}\times{\mathbb{R}}^{d}\times {\mathbb{R}}$\ be
${\cal{P}}\otimes {\cal{B}}({\mathbb{R}}) \otimes
{\cal{B}}({\mathbb{R}}^d) \otimes {\cal{B}}({\mathbb{R}})$\
measurable and satisfy
 \vskip0.2cm

{\rm(i)} There exists a constant $C\ge 0$  such that, P-a.s., for
all $t\in [0, T],\ y_{1}, y_{2}\in {\mathbb{R}},\ z_{1}, z_{2}\in
{\mathbb{R}}^d, \ k_{1}, k_{2}\in {\mathbb{R}},\\ \mbox{ }\hskip2cm
|h(t, y_{1}, z_{1}, k_{1}) - h(t, y_{2}, z_{2}, k_{2})|\leq
C(|y_{1}-y_{2}| + |z_{1}-z_{2}|+|k_{1}-k_{2}|).$
 \vskip0.2cm

{\rm(ii)} $h(\cdot,0,0,0)\in {\cal{H}}^{2}(0,T;{\mathbb{R}})$.
\vskip0.2cm

{\rm(iii)} $k\rightarrow h(t,y,z,k)$\ is non-decreasing, for all
$(t,y,z)\in [0, T]\times {\mathbb{R}}\times{\mathbb{R}}^{d}$.
\vskip0.2cm Furthermore, let $l:\Omega\times [0, T]\times
E\rightarrow {\mathbb{R}}$\ be ${\cal{P}}\otimes {\cal{B}}(E)$\
measurable and satisfy
$$0\leq l_t(e)\leq C(1\wedge|e|),\ \ e\in E.$$
We set
$$
g(t,\omega, y, z, \varphi)=h(t,\omega, y,
z,\int_E\varphi(e)l_t(\omega,e)\lambda(de)),
$$
$\mbox{for}\ (t, \omega, y, z,\varphi)\in [0, T]\times \Omega \times
{\mathbb{R}}\times{\mathbb{R}}^{d}\times L^2(E, {\cal{B}}(E),
\lambda; {\mathbb{R}}).$

Let $\xi',\ \xi\in L^2(\O, {\cal{F}}_T, P)$\ and $g'$\ satisfies (A1)
and (A2).

We denote by $(y,z,k)$\ (resp., $(y',z',k')$) the unique solution of
equation (2.1) with the data $(\xi, g)$\ (resp., $(\xi', g')$ ).
If \vskip0.2cm

{\rm(iv)} $\xi\geq \xi', \ \mbox{a.s.};$
 \vskip0.2cm

{\rm(v)} $g(t,y,z,k)\geq g'(t,y,z,k),  \ \mbox{a.s., a.e., for any}\
(y, z, k)\in {\mathbb{R}}\times{\mathbb{R}}^{d}\times L^2(E,
{\cal{B}}(E), \lambda; {\mathbb{R}}) $. \vskip0.2cm \noindent Then,
we have: $y_t\geq y'_{t},\ a.s.$, for all $t\in [0, T].$\ And if, in
addition, we also assume that $P(\xi_1 > \xi_2)> 0$, then $P\{y_t>
y'_{t}\}>0, \ 0 \leq t \leq T,$\ and in particular, $ y_0> y'_{0}.$
\el

Using the notation introduced in Lemma 2.1 we now suppose that, for
some $g: \Omega\times[0, T]\times{\mathbb{R}}
\times{\mathbb{R}}^{d}\times L^2(E, {\cal{B}}(E), \lambda;
{\mathbb{R}})\longrightarrow {\mathbb{R}}$\ satisfying (A1) and (A2)
and for some $i\in \{1, 2\}$, the drivers $g_i$\ are of the form
$$g_i(s, y_s^i, z_s^i, k_s^i)=g(s, y_s^i, z_s^i, k_s^i)+\varphi_i(s),\ \ \mbox{dsdP-a.e.},$$
where $\varphi_i\in {\cal{H}}^{2}(0,T;{\mathbb{R}}).$\
Then, for terminal values $\xi_1,\ \xi_2\ \mbox{belonging to}\
L^{2}(\Omega, {\cal{F}}_{T}, P)$\ we have the following

 \bl The difference of the solutions $(y^1, z^1, k^1)$ and $(y^2, z^2, k^2)$ of BSDE (2.1) with the data
 $(\xi_1, g_1)$\ and $(\xi_2, g_2)$, respectively, satisfies
 the following estimate:
 $$
  \begin{array}{ll}
  &|y^1_t-y^2_t|^2+\frac{1}{2}E[\int^T_te^{\beta(s-t)}(|
  y^1_s-y^2_s|^2+ |z^1_s-z^2_s|^2)ds|{\cal{F}}_t]\\
  &\ \ \hskip1cm+\frac{1}{2}E[\int^T_t\int_Ee^{\beta(s-t)}|
  k^1_s(e)-k^2_s(e)|^2\lambda(de)ds|{\cal{F}}_t]
     \\
  \leq& E[e^{\beta(T-t)}|\xi_1-\xi_2|^2|{\cal{F}}_t]+ E[\int^T_te^{\beta(s-t)}
           |\varphi_1(s)-\varphi_2(s)|^2ds|{\cal{F}}_t],\ \mbox{P-a.s.,\ for all}\ 0\leq t\leq T,
  \end{array}
  $$
where $\beta\geq 2+2C+4C^2$.
  \el
For the proof the reader is referred to Barles, Buckdahn and
Pardoux~\cite{BBP}, Proposition 2.2.

\section{\large{A DPP for stochastic differential games with jumps}}

\hskip1cm Now we  begin to consider the stochastic differential
games with jumps under our setting.

The set of admissible control processes ${\mathcal{U}}$ (resp.,
${\mathcal{V}}$) for the first (resp., second) player is the set of
all U (resp., V)-valued ${\mathcal{F}}_t$-predictable processes. The
control state spaces U and V are supposed to be compact metric
spaces.

For given admissible controls $u(\cdot)\in {\mathcal{U}}$ and
$v(\cdot)\in {\mathcal{V}}$, the corresponding orbit which regards $t$
as the initial time and $\zeta \in L^2 (\Omega ,{\mathcal{F}}_t,
P;{\mathbb{R}}^n)$ as the initial state is defined by the solution
of the following SDE with jump:
  \be
  \left \{
  \begin{array}{llll}
  dX^{t,\zeta ;u, v}_s & = & b(s,X^{t,\zeta; u,v}_s, u_s, v_s) ds +
          \sigma(s,X^{t,\zeta; u,v}_s, u_s, v_s) dB_s\\
          &&+\int_E\gamma(s,X^{t,\zeta ;u, v}_{s-},u_s,v_s,e)\widetilde{\mu}(ds,de), \hskip1cm  s\in [t,T],\\
   X^{t,\zeta ;u, v}_t  & = & \zeta,
   \end{array}
   \right.
  \ee
where the mappings
  $$
  \begin{array}{llll}
  &   b:[0,T]\times {\mathbb{R}}^n\times U\times V \rightarrow {\mathbb{R}}^n \
  ,\ \ \ \ \ \   \sigma: [0,T]\times {\mathbb{R}}^n\times U\times V\rightarrow {\mathbb{R}}^{n\times d}, \\
  &   \gamma:[0,T]\times {\mathbb{R}}^n\times U\times V\times E \rightarrow
  {\mathbb{R}}^n,
     \end{array}
  $$
  satisfy the following conditions:
  $$
  \begin{array}{ll}
 \rm{(i)}& \mbox{For every fixed}\ (x, e)\in {\mathbb{R}}^n\times E,\ b(.,x,
 .,.), \sigma(.,x,
 .,.)\ \mbox{and}\ \gamma (.,x,
 .,.,e)\ \mbox{are }\\
  &\mbox{continuous in}\ (t,u,v);\\
 \rm{(ii)}&\mbox{There exists a constant }C>0\ \mbox{such that, for all}\ t\in [0,T],\ x, x'\in {\mathbb{R}}^n,\ u \in U,\ v \in V, \\
   &\hskip1cm |b(t,x,u,v)-b(t,x',u ,v)|+ |\sigma(t,x,u,v)-\sigma(t,x',u, v)|\leq C|x-x'|.\\
\mbox{(iii)}&\mbox{There exists }\rho: E \rightarrow
{\mathbb{R}}^{+}\mbox{ with }\int_E \rho^2 (e)\lambda(de)<+\infty,
                              \mbox{ such that,}\\[.1cm]
   & \mbox{for any }t\in [0,T], x, y \in {\mathbb{R}}^n, u\in U, v\in V\mbox{ and }\ e \in E,\\[.1cm]
  & \begin{array}{rcl}
 |\gamma (t,x, u, v, e)- \gamma (t,y, u, v,e)| & \leq & \rho (e)|x-y|, \\[.1cm]
  |\gamma (t,0, u, v,e)| & \leq & \rho(e).
  \end{array}
  \end{array}
  \eqno{\mbox{(H3.1)}}
  $$

From (H3.1) we  get the global linear growth conditions of b and
$\sigma$, i.e., the existence of some $C>0$\ such that, for all $0
\leq t \leq T,\ u\in U,\ v \in V,\  x\in {\mathbb{R}}^n $,
  \be\begin{array}{rcl}
  |b(t,x,u,v)| +|\sigma (t,x,u,v)| &\leq & C(1+|x| );\\
  |\gamma (t,x, u, v,e)| &\leq & \rho(e) (1+|x| ).
  \end{array}
  \ee
Obviously, under the above assumptions, for any $u(\cdot)\in
{\mathcal{U}}$ and $v(\cdot)\in {\mathcal{V}}$, SDE (3.1) has a
unique strong solution. Moreover, there exists $C\in \mathbb{R}^+$\
such that, for any $t \in [0,T]$, $u(\cdot)\in {\mathcal{U}},
v(\cdot)\in {\mathcal{V}}$\ and $ \zeta, \zeta'\in L^2 (\Omega
,{\mathcal{F}}_t,P;{\mathbb{R}}^n),$\
we have the following estimates, P-a.s.:\\
 \be
\begin{array}{rcl}
E[\sup \limits_{s\in [t,T]}|X^{t,\zeta; u, v}_s -X^{t,\zeta';u,
v}_s|^2|{{\mathcal{F}}_t}]
& \leq & C|\zeta -\zeta'|^2, \\
E[ \sup \limits_{s\in [t,T]} |X^{t,\zeta
;u,v}_s|^2|{{\mathcal{F}}_t}] & \leq &
                        C(1+|\zeta|^2).
\end{array}
\ee The constant $C$ depends only on the Lipschitz
 and the linear growth constants of $b,\ \sigma$\ and $\gamma$
with respect to $x$.

Let now be given three measurable functions
$$
\Phi: {\mathbb{R}}^n \rightarrow {\mathbb{R}},\  f:[0,T]\times
{\mathbb{R}}^n \times {\mathbb{R}} \times {\mathbb{R}}^d \times
{\mathbb{R}}\times U \times V \rightarrow {\mathbb{R}},\ l:
{\mathbb{R}}^n\times E\rightarrow {\mathbb{R}}
$$
which satisfy the following conditions:
$$
\begin{array}{ll}
\rm{(i)}& \mbox{For every fixed}\ (x, y, z, k)\in {\mathbb{R}}^n
\times {\mathbb{R}} \times {\mathbb{R}}^d\times {\mathbb{R}},\ f(.,
x, y, z, k,.,.)\
\mbox{is continuous in}\ \\
&(t,u,v)\ \mbox{and there exists a constant}\ C>0 \ \mbox{such that,
for all}\ t\in [0,T],\ x, x'\in {\mathbb{R}}^n,\\\
& y, y'\in {\mathbb{R}},\ z, z'\in {\mathbb{R}}^d,\ k, k'\in {\mathbb{R}}, u \in U \ \mbox{and}\ v \in V,\\
&\hskip2cm\begin{array}{l}
|f(t,x,y,z,k,u,v)-f(t,x',y',z',k',u,v)| \\
\hskip4cm \leq C(|x-x'|+|y-y'| +|z-z'|+|k-k'|);
\end{array}\\
\rm{(ii)}&k\rightarrow f(t,x,y,z,k,u,v)\ \mbox{is non-decreasing},\
\mbox{for all}\ (t,x,y,z,u,v)\in [0, T]\times {\mathbb{R}}^n\times
{\mathbb{R}}\\
&\times{\mathbb{R}}^{d}\times U\times V;\\
\rm{(iii)}& \mbox{There exists a constant}\ C>0\ \mbox{such that},
\\
&\hskip2cm\begin{array}{rcl} 0\leq
l(x,e)&\leq& C(1\wedge|e|),\ \ x\in {\mathbb{R}}^n,\ e\in E,\\
|l(x,e)-l(x',e)|&\leq &C|x-x'|(1\wedge|e|),\ x, x'\in
{\mathbb{R}}^n,\ e\in E;
\end{array}\\
\rm{(iv)}&\mbox{There exists a constant}\ C>0 \ \mbox{such that, for
all}\ x, x'\in {\mathbb{R}}^n,\\
 &\mbox{  }\hskip3cm |\Phi (x) -\Phi (x')|\leq C|x-x'|.
 \end{array}
 \eqno {\mbox{(H3.2)}}
 $$
From (H3.2) we see that $f$\ and $\Phi$\ also satisfy the global
linear growth condition in $x$, i.e., there exists some $C>0$\ such
that, for all $0 \leq t \leq T,\  u\in U,\ v \in V,\ x\in
{\mathbb{R}}^n $,
  \be
   |f(t,x,0,0,0,u,v)|+|\Phi (x)| \leq C(1+|x|).
   \ee
For any $u(\cdot) \in {\mathcal{U}}, $\ $v(\cdot) \in
{\mathcal{V}}$\ and $\zeta \in L^2
(\Omega,{\mathcal{F}}_t,P;{\mathbb{R}}^n)$,  the mappings $\xi:=
\Phi(X^{t,\zeta; u, v}_T)$ and $g(s,y,z,k):= f(s,X^{t,\zeta; u,
v}_s,y,z,\int_Ek(e)l(X_s^{t, \zeta;u, v}, e)\lambda(de),u_s,v_s),\ \
(s, y, z, k)\in [0, T]\times {\mathbb{R}}\times {\mathbb{R}}^d\times
L^2(E, {\cal{B}}(E), \lambda;{\mathbb{R}})$ satisfy the conditions
of Lemma 2.1 on the interval $[t, T]$. Therefore, there exists a
unique solution to the following BSDE:
      \be
   \left \{\begin{array}{rcl}
   -dY^{t,\zeta; u, v}_s & = & f(s,X^{t,\zeta; u, v}_s, Y^{t,\zeta; u, v}_s, Z^{t,\zeta; u,
   v}_s,\int_EK^{t,\zeta; u, v}_s(e)l(X^{t,\zeta; u,
   v}_s,e)\lambda(de),u_s, v_s) ds\\
   &&-Z^{t,\zeta; u, v}_s dB_s-\int_EK^{t,\zeta; u, v}_s(e)\widetilde{\mu}(ds,de),\\
      Y^{t,\zeta; u, v}_T  & = & \Phi (X^{t,\zeta; u, v}_T),
   \end{array}\right.
   \ee
where $X^{t,\zeta; u, v}$\ is introduced by equation (3.1).

Note that in (3.5) and in the sequel, $f$\ depends on $K$\ in a very
specific way in order to make full use of the comparison
theorem-Lemma 2.2.

 Moreover, in analogy to Proposition 6.1 in the Appendix, we can see that
there exists some constant $C>0$\ such that, for all $0 \leq t
\leq T,\ \zeta, \zeta' \in L^2(\Omega ,
{\mathcal{F}}_t,P;{\mathbb{R}}^n),\ u(\cdot) \in {\mathcal{U}}\
\mbox{and}\ v(\cdot) \in {\mathcal{V}},$\ P-a.s.,
 \be
\begin{array}{ll}
 {\rm(i)} & |Y^{t,\zeta; u, v}_t -Y^{t,\zeta'; u, v}_t| \leq C|\zeta -\zeta'|; \\
 {\rm(ii)} & |Y^{t,\zeta; u, v}_t| \leq C (1+|\zeta|). \\
\end{array}
\ee

Now, similar to~\cite{BL} and~\cite{FS}, we introduce the following
subspaces of admissible controls and the definition of admissible
strategies for the game:

\noindent\bde\ An admissible control process $u=\{u_r, r\in [t,
s]\}$ (resp., $v=\{v_r, r\in [t, s]\}$) for Player I (resp., II) on
$[t, s] (t<s\leq T)$\ is an ${\mathcal{F}}_r$-predictable process
taking values in U (resp., V). The set of all admissible controls
for Player I (resp., II) on $[t, s]$ is denoted by\
${\mathcal{U}}_{t, s}$\ (resp., ${\mathcal{V}}_{t, s}).$\ We
identify both processes $u$\ and $\bar{u}$\ in\ ${\mathcal{U}}_{t,
s}$\ and write $u\equiv \bar{u}\ \mbox{on}\ [t, s],$\ if
$P\{u=\bar{u}\ \mbox{a.e. in}\ [t, s]\}=1.$\ Similarly we interpret
$v\equiv \bar{v}\ \mbox{on}\ [t, s]$\ in ${\mathcal{V}}_{t, s}$.
\ede
 Finally, we still have to define the admissible strategies for the
game.

\bde A nonanticipative strategy for Player I on $[t, s] (t<s\leq
T)$ is a mapping $\alpha: {\mathcal{V}}_{t, s}\longrightarrow
{\mathcal{U}}_{t, s}$ such that, for any
${\mathcal{F}}_r$-stopping time $S: \Omega\rightarrow [t, s]$\ and
any $ v_1, v_2 \in {\mathcal{V}}_{t, s}$\ with $ v_1\equiv v_2\
\mbox {on}\ \textbf{[\![}t, S\textbf{]\!]},$ it holds
$\alpha(v_1)\equiv \alpha(v_2)\ \mbox {on}\ \textbf{[\![}t,
S\textbf{]\!]}$.\ Nonanticipative strategies for Player II on $[t,
s]$, $\beta: {\mathcal{U}}_{t, s}\longrightarrow {\mathcal{V}}_{t,
s}$,  are defined similarly. The set of all nonanticipative
strategies $\alpha: {\mathcal{V}}_{t,s}\longrightarrow
{\mathcal{U}}_{t,s}$ for Player I on $[t, s]$ is denoted by
${\cal{A}}_{t,s}$. The set of all nonanticipative strategies
$\beta: {\mathcal{U}}_{t,s}\longrightarrow {\mathcal{V}}_{t,s}$
for Player II on $[t, s]$ is denoted by ${\cal{B}}_{t,s}$.
\\
(\mbox{Recall that}\ $\textbf{[\![}t,
S\textbf{]\!]}=\{(r,\omega)\in [0, T]\times \Omega, t\leq r\leq
S(\omega)\})$.\ede

Given the control processes $u(\cdot)\in {\mathcal{U}}_{t,T}$\ and $
v(\cdot)\in {\mathcal{V}}_{t,T} $\ we introduce the following
associated cost functional
 \be
J(t, x; u, v):= Y^{t, x; u, v}_t,\ (t, x)\in [0, T]\times
{\mathbb{R}}^n,\ee where the process $Y^{t, x; u, v}$ is defined by
BSDE (3.5).

 \noindent Similarly to the
proof of Theorem 6.1 in the Appendix, we can get that, for any
$t\in[0, T]$\ and $\zeta \in L^2 (\Omega ,{\mathcal{F}}_t ,P;
{\mathbb{R}}^n)$,
 \be J(t, \zeta; u, v) = Y^{t,\zeta; u, v}_t,\
 \mbox{P-a.s.}
\ee
 Being particularly interested in the case of a deterministic $\zeta$, i.e., $\zeta=x\in {\mathbb{R}}^n$,
 we define the lower value function of our stochastic
differential game \be W(t,x):= \mbox{essinf}_{\beta \in
{\cal{B}}_{t,T}}\mbox{esssup}_{u \in {\mathcal{U}}_{t,T}}J(t,x;
u,\beta(u)) \ee
 and its upper value function
  \be U(t,x):= \mbox{esssup}_{\alpha \in
{\cal{A}}_{t,T}}\mbox{essinf}_{v \in {\mathcal{V}}_{t,T}}J(t,x;
\alpha(v),v). \ee

\br (1) Here the essential infimum and the essential supremum should
be understood as one with respect to indexed families of random
variables (see, e.g., Dunford and Schwartz~\cite{DS},
Dellacherie~\cite{D} or the Appendix in Karatzas and
Shreve~\cite{KS2} for detailed discussions). The reader is also
referred to Remark 3.1 in~\cite{BL}.\\
(2) Let us point out that under our conditions (H3.1)-(H3.2) the lower value
function $W(t,x)$\ and the upper value function $U(t,x)$ are well defined and, a priorily, bounded, ${\cal F}_t$-measurable random variables. However, we show below that they are indeed deterministic functions. Such a result was already got in the case of stochastic differential games only driven by a Brownian motion (see~\cite{BL}). However, here, in presence of an additional driving compensated Poisson random measure, the argument of the Girsanov transformation employed in~\cite{BL} doesn't work anymore and has to be replaced by a quite different transformation argument. In what follows we concentrate on the study of $W$, the upper value function $U$ can be investigated in a similar manner.\er
 \bp For any $(t, x)\in [0, T]\times {\mathbb{R}}^n$,
we have $W(t,x)=E[W(t,x)]$, P-a.s. Thus, let $W(t,x)$\ identify with
its deterministic version $E[W(t,x)],$\ $W:[0,
T]\times {\mathbb{R}}^n\longrightarrow {\mathbb{R}}$ is a
deterministic function.\ep

The proof will be split into two lemmas.

\bl
Let $(t,x)\in [0, T]\times {\mathbb{R}}^n$\ and $\tau: \Omega\rightarrow\Omega$ be an invertible ${\cal{F}}$-${\cal{F}}$ measurable transformation such that
$$
\begin{array}{lll}
& \rm{i)}\ \tau\  \mbox{and}\  \tau^{-1}: \Omega\rightarrow\Omega\  \mbox{are}\  {\cal{F}}_t-{\cal{F}}_t\  \mbox{measurable};\\
& \rm{ii)}\ (B_s-B_t)\circ\tau=B_s-B_t,\ s\in[t, T];\\
&\ \ \ \ \mu((t, s]\times A)\circ\tau=\mu((t, s]\times A),\ s\in[t, T],\ A\in {\cal{B}}(E);\\
& \rm{iii)}\ \mbox{the law}\ P\circ [\tau]^{-1}\ \mbox{of}\ \ \tau\ \mbox{is equivalent to the underlying probability measure}\ P.
\end{array}
$$
Then, $W(t,x)\circ\tau=W(t,x),$ P-a.s.
\el
\noindent \textbf{Proof}: We split now
the proof in the following steps:
 \vskip0.1cm
\noindent $1^{st}$ step: For any $u\in {\mathcal{U}}_{t,T}, \ v\in
{\mathcal{V}}_{t,T}, \ J(t, x; u,v)\circ\tau= J(t, x;
u(\tau),v(\tau)),\ \mbox{P-a.s.}$ \vskip0.1cm
 Indeed, we apply the transformation $\tau$ to SDE (3.1) (with
 $\zeta=x$) and compare the obtained equation with the SDE
 obtained from (3.1) by substituting the controlled processes $u(\tau),\ v(\tau)$\ for $u$\ and $v$. Then, from the uniqueness of the solution of (3.1) we get $X_s^{t,x; u,v}(\tau)=X_s^{t,x;
u(\tau),v(\tau)},$ $ \mbox{for any}\ s\in [t, T],\
\mbox{P-a.s.}$\ Furthermore, by a similar transformation
argument we obtain from the uniqueness of the solution of BSDE (3.5),
$$Y_s^{t,x; u,v}(\tau)=Y_s^{t,x; u(\tau),v(\tau)},\ \mbox{for
any}\ s\in [t, T],\ \mbox{P-a.s.,}$$
$$Z_s^{t,x; u,v}(\tau)=Z_s^{t,x; u(\tau),v(\tau)},\  \mbox{dsdP-a.e. on}\ [t, T]\times\Omega,$$
$$K_s^{t,x; u,v}(\tau)=K_s^{t,x; u(\tau),v(\tau)},\  ds\lambda(de)\mbox{dP-a.e. on}\ [t, T]\times E\times\Omega.$$
Consequently, in particular, we have $$J(t, x; u,v)(\tau)= J(t, x;
u(\tau),v(\tau)),\ \mbox{P-a.s.}$$
 \vskip0.1cm
\noindent $2^{nd}$ step: For $\beta\in {\cal{B}}_{t,T},$\ let $\widehat{\beta}(u):=\beta(u(\tau^{-1}))(\tau),\ u\in
{\mathcal{U}}_{t,T}.$\ Then, $\widehat{\beta}\in {\cal{B}}_{t,T}.$
\vskip0.1cm Obviously, $\widehat{\beta}$\ maps ${\mathcal{U}}_{t,T}$\ into
${\mathcal{V}}_{t,T}$.\ Moreover, this mapping $\widehat{\beta}$ is nonanticipating.
Indeed, let $S: \Omega\rightarrow [t, T]$\ be an
${\mathbb{F}}$-stopping time and $ u_1, u_2 \in
{\mathcal{U}}_{t, T}$\ such that $ u_1\equiv u_2\ \mbox {on}\
\textbf{[\![}t,\ S\textbf{]\!]}.$\ Then, obviously, $
u_1(\tau^{-1})\equiv u_2(\tau^{-1})\ \mbox {on}\ \textbf{[\![}t,\
S(\tau^{-1})\textbf{]\!]}$ (notice that $ S(\tau^{-1})\ \mbox{is
still an}\ {\mathbb{F}}$-stopping time. For this we use that the assumptions i) and ii) imply that $\tau({\cal{F}}_s):=\{\tau(A),\ A\in {\cal{F}}_s\}={\cal{F}}_s,\ s\in[t, T]$ ). Thus, because $\beta\in {\cal{B}}_{t,T}$,\ we have $\beta(u_1(\tau^{-1}))\equiv \beta(u_2(\tau^{-1}))\ $ $ \mbox
{on}\ \textbf{[\![}t,\ S(\tau^{-1})\textbf{]\!]}$. Therefore,
$$\widehat{\beta}(u_1)=\beta(u_1(\tau^{-1}))(\tau)\equiv \beta(u_2(\tau^{-1}))(\tau)=\widehat{\beta}(u_2)\ \mbox
{on}\ \textbf{[\![}t,\ S\textbf{]\!]}.$$
\vskip0.1cm
\noindent$3^{rd}$ step: For all $\beta\in
{\mathcal{B}}_{t, T}$\ we have:
$$(\mbox{esssup}_{u \in {\mathcal{U}}_{t,T}}J(t,x;
u,\beta(u)))(\tau)=\mbox{esssup}_{u \in
{\mathcal{U}}_{t,T}}(J(t,x; u,\beta(u))(\tau)),\ \mbox{P-a.s.}.
$$

Indeed, with the notation $I(t,x;\beta):=\mbox{esssup}_{u \in
{\mathcal{U}}_{t,T}}J(t,x; u,\beta(u)),\ \beta\in {\mathcal{B}}_{t,
T},$\ we have \ $I(t,x;\beta)\geq J(t,x; u,\beta(u)),$\ and thus
$I(t,x;\beta)(\tau)\geq J(t,x; u,\beta(u))(\tau), \
\mbox{P-a.s.,\ for}$ $\mbox{ all}\ u\in {\mathcal{U}}_{t,T}$\
(recall that $P\circ\tau^{-1}$\ is equivalent to $P$\ due to assumption iii)). On the
other hand, for any random variable $\zeta$\ satisfying $\zeta\geq
J(t,x; u,\beta(u))(\tau),$\ and hence also $\zeta(\tau^{-1})\geq
J(t,x; u,\beta(u)), \ \mbox{P-a.s.,\ for}\ \mbox{ all}\ u\in
{\mathcal{U}}_{t,T},$\ we have\ $\zeta(\tau^{-1})\geq I(t,x;\beta),
\ $ $ \mbox{P-a.s.,}$\ i.e., $\zeta\geq I(t,x;\beta)(\tau), \
\mbox{P-a.s.}$\ Consequently,
$$I(t,x;\beta)(\tau)=\mbox{esssup}_{u \in
{\mathcal{U}}_{t,T}}(J(t,x; u,\beta(u))(\tau)),\ \mbox{P-a.s.}$$
 \vskip0.1cm

 \noindent$4^{th}$ step: $W(t,x)$\ is invariant with respect
 to the transformation $\tau$, i.e.,
  $$W(t,x)(\tau)=W(t,x), \ \mbox{P-a.s.}$$

Indeed, similarly to the third step we can show that:
$$(\mbox{essinf}_{\beta \in
{\mathcal{B}}_{t,T}}I(t,x;\beta))(\tau)=\mbox{essinf}_{\beta \in
{\mathcal{B}}_{t,T}}(I(t,x; \beta)(\tau)),\ \mbox{P-a.s.}
$$
Then, from the first step to the third step we have,
 $$
   \begin{array}{rcl}
   W(t,x)(\tau) & = & \mbox{essinf}_{\beta \in
{\mathcal{B}}_{t,T}}\mbox{esssup}_{u \in
{\mathcal{U}}_{t,T}}(J(t,x; u,\beta(u))(\tau))\\
       & = &  \mbox{essinf}_{\beta \in
{\mathcal{B}}_{t,T}}\mbox{esssup}_{u \in
{\mathcal{U}}_{t,T}}J(t,x; u(\tau),\widehat{\beta}(u(\tau)))\\
& = &  \mbox{essinf}_{\beta \in
{\mathcal{B}}_{t,T}}\mbox{esssup}_{u \in
{\mathcal{U}}_{t,T}}J(t,x; u,\widehat{\beta}(u))\\
& = &  \mbox{essinf}_{\beta \in
{\mathcal{B}}_{t,T}}\mbox{esssup}_{u \in
{\mathcal{U}}_{t,T}}J(t,x; u,\beta(u))\\
& = &W(t,x),\ \mbox{P-a.s.,}
   \end{array}
$$
where we have used
$\{u(\tau)\ |\ u(\cdot) \in {\mathcal{U}}_{t,T}\}={\mathcal{U}}_{t,T},\
\{\widehat{\beta}\ |\ \beta \in {\mathcal{B}}_{t,T} \}={\mathcal{B}}_{t,T}$\
in order to obtain the both latter equalities. \endpf

\vskip0.1cm

Now let $ \ell\geq 1$. We define the transformation $\tau'_\ell: \Omega_1\rightarrow\Omega_1$\ such that, for all $\omega_1\in \Omega_1=C_0({\mathbb{R}};{\mathbb{R}}^d)$,
\be
 \begin{array}{lll}
 &(\tau'_\ell\omega_1)((t-\ell, r])=\omega_1((t-2\ell, r-\ell])(:=\omega_1(r-\ell)-\omega_1(t-2\ell));\\
 &(\tau'_\ell\omega_1)((t-2\ell, r-\ell])=\omega_1((t-\ell, r]),\ \mbox{for}\ r\in [t-\ell, t];\\
 &(\tau'_\ell\omega_1)((s, r])=\omega_1((s, r]),\ (s, r]\cap (t-2\ell, t]=\emptyset;\\
 &(\tau'_\ell\omega_1)(0)=0.
 \end{array}
\ee
Moreover, for $p\in \Omega_2,\ p=\Sigma_{x\in D_p}p(x)\delta_x$, we put:
 $$\tau''_\ell p:=\Sigma_{x\in D_p\cap (t-2\ell, t]^c}p(x)\delta_x+\Sigma_{x\in D_p\cap (t-\ell, t]}p(x)\delta_{x-\ell}+\Sigma_{x\in D_p\cap (t-2\ell, t-\ell]}p(x)\delta_{x+\ell}.$$

\noindent Obviously, $\tau''_\ell: \Omega_2\rightarrow \Omega_2$\ is a bijection, ${\tau''_\ell}^{-1}=\tau''_\ell$, which preserves the measure $P_2\circ [\tau''_\ell]^{-1}=P_2$. Moreover,
 \be
 \begin{array}{lll}
 &\mu(\tau''_\ell p ;\ (t-\ell, r]\times \Delta)=\mu(p ;\ (t-2\ell, r-\ell]\times \Delta), \ r\in (t-\ell, t],\ \Delta\in {\cal B}(E);\\
 &\mu(\tau''_\ell p ;\ (t-2\ell, r-\ell]\times \Delta)=\mu(p ;\ (t-\ell, r]\times \Delta),\ r\in (t-\ell, t],\ \Delta\in {\cal B}(E);\\
 &\mu(\tau''_\ell p ;\ (s, r]\times \Delta)=\mu(p ;\ (s, r]\times \Delta),\ (s, r]\cap (t-2\ell, t]=\emptyset, \ \Delta\in {\cal B}(E).\\
 \end{array}
\ee
Thus, the transformation $\tau_\ell: \Omega\rightarrow\Omega$, $\tau_\ell\omega:=(\tau'_\ell\omega_1,\ \tau''_\ell p),\ \omega=(\omega_1,\ p)\in \Omega=\Omega_1\times\Omega_2$, satisfies the assumptions i), ii), iii) of Lemma 3.1. Therefore, $W(t, x)(\tau_\ell)=W(t, x),\ \mbox{P-a.s.},\ \ell \geq 1.$\ The proof of Proposition 3.1 will be completed by the following auxiliary Lemma 3.2.

\bl Let $\zeta \in L^\infty(\Omega, {\cal F}_t, P)$\ be such that, for all $\ell \geq 1$\ natural number, $\zeta(\tau_\ell)=\zeta$, P-a.e. Then, there exists some real $C$\ such that $\zeta=C$, P-a.s.
\el

\noindent \textbf{Proof}: To simplify the notation we introduce the Brownian motion $B'_r:=B_t-B_{t-r},\ r\geq 0,$\ and the Poisson random measure
 $$\upsilon([0, r]\times \Delta):=\mu([t-r, t]\times \Delta),\ \  r\geq 0,\ \Delta\in {\cal B}(E),\ $$
 $\mbox{on}\ {\mathbb{R}}_{+}\times E.$\ Then,
 $$
 \begin{array}{lll}
 &B'_r(\tau'_\ell\omega_1)=B'_{\ell+r}(\omega_1)-B'_{\ell}(\omega_1);\\
 &(B'_{\ell+r}-B'_{\ell})(\tau'_\ell\omega_1)=B'_r(\omega_1),\ r\in [0, \ell];\\
 &B'_{r}(\tau'_\ell\omega_1)-B'_{s}(\tau'_\ell\omega_1)=B'_{r}(\omega_1)-B'_{s}(\omega_1),\ (s, r]\cap [0, 2\ell]=\emptyset; \\
 \end{array}
 $$
 and
 $$
 \begin{array}{lll}
 &\upsilon(\tau''_\ell p ;\ [0, r]\times \Delta)=\upsilon(p ;\ [\ell, \ell+r]\times \Delta);\\
 &\upsilon(\tau''_\ell p ;\ [\ell, \ell+r]\times \Delta)=\upsilon(p ;\ [0, r]\times \Delta),\ r\in [0, \ell],\ \Delta\in {\cal B}(E);\\
 &\upsilon(\tau''_\ell p ;\ [s, r]\times \Delta)=\upsilon(p ;\ [s, r]\times \Delta),\ [s, r]\cap [0, 2\ell]=\emptyset, \ 0\leq s \leq r,\ \Delta\in {\cal B}(E).\\
 \end{array}
 $$
 Moreover, we put ${\cal F}'_{s, r}:=\sigma\{B'_{s'}-B'_{s},\ \upsilon([s, s']\times \Delta),\ s'\in [s, r],\ \Delta\in {\cal B}(E)\}\vee {\cal N}_P,$
 $0\leq s \leq r<+\infty$. Let $\zeta \in L^\infty(\Omega, {\cal F}_t, P)$\ be such that $\zeta(\tau_\ell)=\zeta$, P-a.s., for all $\ell \geq 1$\ natural number. To prove Lemma 3.2 it suffices to show that $\zeta=E[\zeta]$, P-a.s.

 To this end we consider a random variable of the form $\theta\eta\varsigma$, where $\theta\in L^\infty (\Omega, {\cal F}'_{0, \ell}, P), \ \eta\in L^\infty (\Omega, {\cal F}'_{\ell, 2\ell}, P),$\ and $\varsigma \in L^\infty (\Omega, {\cal F}'_{2\ell, \infty}, P)$. Then, $\theta(\tau_\ell)\in L^\infty (\Omega, {\cal F}'_{\ell, 2\ell}, P),$\ and the random variable $\eta(\tau_\ell)\in L^\infty (\Omega, {\cal F}'_{0, \ell}, P)$\ is independent of ${\cal F}'_{\ell, \infty}$. Consequently, taking into account that $\varsigma(\tau_\ell)=\varsigma$, P-a.e., we have
 $$E[(\theta\eta\varsigma)(\tau_\ell)|{\cal F}'_{\ell, \infty}]=E[\eta(\tau_\ell)]\theta(\tau_\ell)\varsigma=E[\eta]\theta(\tau_\ell)\varsigma=E[\theta\eta\varsigma|{\cal F}'_{0, \ell}\vee {\cal F}'_{2\ell, \infty} ](\tau_\ell),\ \mbox{P-a.s.,}$$
and from the monotone class theorem we conclude that, for all $\theta\in L^1(\Omega, {\cal F}'_\infty, P),$
$$E[\theta(\tau_\ell)|{\cal F}'_{\ell, \infty}]=E[\theta|{\cal F}'_{0, \ell}\vee {\cal F}'_{2\ell, \infty}](\tau_\ell),\ \mbox{P-a.s.,}\ \ell\geq 1.$$
Thus, $E[\zeta|{\cal F}'_{\ell, \infty}]=E[\zeta(\tau_\ell)|{\cal F}'_{\ell, \infty}]=E[\zeta|{\cal F}'_{0, \ell}\vee {\cal F}'_{2\ell, \infty}](\tau_\ell),\ \mbox{P-a.s.},\ \ell\geq 1.$\

As $E[\zeta|{\cal F}'_{\ell, \infty}]\rightarrow E[\zeta|\cap_{\ell\geq 1}{\cal F}'_{\ell, \infty}]=E[\zeta],$\ as $\ell\rightarrow \infty$, P-a.s., and in $L^1$, it follows that
$$
\begin{array}{lll}
&\ \ E[|E[\zeta|{\cal F}'_{0, \ell}]-E[\zeta]|]\\
&=E[|E[E[\zeta|{\cal F}'_{0, \ell}\vee {\cal F}'_{2\ell, \infty}]-E[\zeta]|{\cal F}'_{0, \ell}]|]\\
&\leq E[|E[\zeta|{\cal F}'_{0, \ell}\vee {\cal F}'_{2\ell, \infty}]-E[\zeta]|]\\
&=E[|E[\zeta|{\cal F}'_{0, \ell}\vee {\cal F}'_{2\ell, \infty}](\tau_\ell)-E[\zeta]|] \ \ \ \ (\mbox{Recall}: P\circ [\tau_\ell]^{-1}=P)\\
&=E[|E[\zeta|{\cal F}'_{\ell, \infty}]-E[\zeta]|]\rightarrow 0,\ \mbox{as}\ \ell\rightarrow +\infty.
\end{array}
$$
Consequently, $E[\zeta|{\cal F}'_{0, \ell}]\rightarrow E[\zeta]$\ in $L^1$, as $\ell\rightarrow \infty$. But, on the other hand, $E[\zeta|{\cal F}'_{0, \ell}]\rightarrow E[\zeta|{\cal F}'_{0, \infty}]=\zeta,$\ in $L^1$, as $\ell\rightarrow \infty$. This shows that $\zeta=E[\zeta]$, P-a.s. \endpf

\br
From Lemma 3.2 we know $W(t,x)$\ is independent of
 ${\mathcal{F}}_{T}.$
\er
  \vskip0.1cm

The first property of the lower value function $W(t,x)$ which we
present is an immediate consequence of its definition (3.9) and the
Lipschitz property (3.6) of the cost functionals.

\bl\mbox{  }There exists a constant $C>0$\ such that, for all $ 0
\leq t \leq T,\ x, x'\in {\mathbb{R}}^n$,\be
\begin{array}{llll}
&{\rm(i)} & |W(t,x)-W(t,x')| \leq C|x-x'|;  \\
&{\rm(ii)} & |W(t,x)| \leq C(1+|x|).
\end{array}
\ee \el \endpf \vskip0.3cm We now discuss (the generalized) dynamic
programming principle (DPP) for our stochastic differential game
(3.1), (3.5) and (3.9). For this end we have to
 define the family of (backward) semigroups associated with BSDE
 (3.5). This notion of the stochastic backward semigroup was first
 introduced by Peng~\cite{P1} which was applied to study the DPP for
 stochastic control problems in the framework of Brownian motion. Our approach adapts Peng's ideas to the framework
  of stochastic differential games with jumps.

 Given the initial data $(t,x)$, a positive number $\delta\leq T-t$, admissible control
 processes $u(\cdot) \in {\mathcal{U}}_{t, t+\delta},\ v(\cdot) \in {\mathcal{V}}_{t, t+\delta}$\ and a real-valued
 random variable $\eta \in L^2 (\Omega,
{\mathcal{F}}_{t+\delta},P;{\mathbb{R}})$, we put \be G^{t, x; u,
v}_{s,t+\delta} [\eta]:= \tilde{Y}_s^{t,x; u, v},\ \hskip0.5cm
s\in[t, t+\delta], \ee where $(\tilde{Y}_s^{t,x;u, v},
\tilde{Z}_s^{t,x;u, v}, \widetilde{K}^{t,x; u, v}_s )_{t\leq s \leq
t+\delta}$ is the solution of the following BSDE with the time
horizon $t+\delta$:
\be
\left \{\begin{array}{rcl}
 -d\tilde{Y}_s^{t,x;u, v} \!\!\!& = &\!\!\! f(s,X^{t,x;u, v}_s ,\tilde{Y}_s^{t,x;u, v}, \tilde{Z}_s^{t,x;u, v},
            \int_E\widetilde{K}^{t,x; u, v}_s(e)l(X^{t,x; u,v}_s,e)\lambda(de), u_{s}, v_{s})ds \\
  && -\tilde{Z}_s^{t,x; u, v} dB_s-\int_E\widetilde{K}^{t,x; u, v}_s(e)\widetilde{\mu}(ds,de), \hskip 1cm s\in [t,t+\delta],\\
 \tilde{Y}_{t+\delta}^{t,x; u, v}\!\!\! & =& \!\!\!\eta ,
\end{array}\right.
\ee where $X^{t,x;u, v}$\ is the solution of SDE (3.1).

\br When $f$\ is independent of $(y, z, k)$\ it holds that
$$G^{t,x;u,v}_{s,t+\delta}[\eta]=E[\eta + \int_s^{t+\delta}
f(r,X^{t,x;u, v}_r,u_{r}, v_{r})dr|{\cal{F}}_s],\ \ s\in [t,
t+\delta].$$ \er

 Obviously, for the solution $(Y^{t,x;u, v}, Z^{t,x;u, v},$
$ K^{t,x;u, v})$\ of BSDE (3.5) we have \be G^{t,x;u, v}_{t,T} [\Phi
(X^{t,x; u, v}_T)] =G^{t,x;u, v}_{t,t+\delta} [Y^{t,x;u,
v}_{t+\delta}]. \ee Moreover,
$$
\begin{array}{rcl}
 J(t,x;u, v)& = &Y_t^{t,x;u, v}=G^{t,x;u, v}_{t,T} [\Phi (X^{t,x; u, v}_T)]
  =G^{t,x;u,v}_{t,t+\delta} [Y^{t,x;u, v}_{t+\delta}]\\
  &=&G^{t,x;u,v}_{t,t+\delta} [J(t+\delta,X^{t,x;u, v}_{t+\delta};u, v)].
\end{array}
$$

 \bt\mbox{}Under the
assumptions (H3.1) and (H3.2), the lower value function $W(t,x)$
obeys the following DPP : For any $0\leq t<t+\delta \leq T,\ x\in {\mathbb{R}}^n,$
 \be
W(t,x) =\mbox{essinf}_{\beta \in {\mathcal{B}}_{t,
t+\delta}}\mbox{esssup}_{u \in {\mathcal{U}}_{t,
t+\delta}}G^{t,x;u,\beta(u)}_{t,t+\delta} [W(t+\delta,
X^{t,x;u,\beta(u)}_{t+\delta})].
 \ee
  \et

The proof is given in  Section 6.2 of the Appendix since it is quite lengthy.

In Lemma 3.3 we have already seen that the lower value function
$W(t,x)$\ is Lipschitz continuous in $x$, uniformly in $t$. With the
help of Theorem 3.1 we can now also study the continuity property of
$W(t,x)$\ in $t$.
 \bt\mbox{ }Let us suppose that the assumptions (H3.1) and (H3.2)
hold. Then the lower value function $W(t,x)$ is
 $\frac{1}{2}-$H\"{o}lder continuous in $t$: there exists a constant C such that,
  for every $x\in {\mathbb{R}}^n,\ t, t'\in [0, T]$,
  $$
|W(t, x)-W(t', x)|\leq C(1+|x|)|t-t'|^{\frac{1}{2}}.
  $$
  \et
\noindent \textbf{Proof}:  Let $(t, x)\in [0,T]\times
{\mathbb{R}}^n$\ and $\delta>0$\ be arbitrarily given such that
$0<\delta\leq T-t$. Our objective is to prove the following
inequality by using (6.21) and (6.22) in the Appendix:
 \be -C(1+|x|)\delta^{\frac{1}{2}}\leq W(t,x)-W(t+\delta ,x)\leq
C(1+|x|)\delta^{\frac{1}{2}}.
 \ee
 From it we obtain immediately that $W$ is $\frac{1}{2}-$H\"{o}lder continuous in
 $t$. We will only check the second inequality in (3.18), the
 first one can be shown in a similar way. To this end we note that
 due to (6.21), for an arbitrarily small $\varepsilon>0,$
\be W(t,x)-W(t+\delta ,x) \leq I^1_\delta +I^2_\delta
+\varepsilon, \ee
 where
$$
\begin{array}{lll}
I^1_\delta & := & G^{t,x;
u^{\varepsilon},\beta(u^{\varepsilon})}_{t,t+\delta}[W(t+\delta,
X^{t,x; u^{\varepsilon},\beta(u^{\varepsilon})}_{t+\delta})]
                   -G^{t,x;u^{\varepsilon},\beta(u^{\varepsilon})}_{t,t+\delta} [W(t+\delta,x)], \\
I^2_\delta & := & G^{t,x;
u^{\varepsilon},\beta(u^{\varepsilon})}_{t,t+\delta}
[W(t+\delta,x)] -W(t+\delta ,x),
\end{array}
$$
for arbitrarily chosen $\beta\in {\cal{B}}_{t, t+\delta}$\ and
$u^{\varepsilon} \in {\cal{U}}_{t, t+\delta}$\ such that (6.21)
holds. From Lemma 2.3 and the estimate (3.13)-(i) we obtain that, for
some constant $C$ independent of the controls $u^{\varepsilon}\
\mbox{and}\ \ \beta(u^{\varepsilon})$,
$$
\begin{array}{rcl}
|I^1_\delta | &\leq& [CE(|W(t+\delta ,X^{t,x;
u^{\varepsilon},\beta(u^{\varepsilon})}_{t+\delta})
                  -W(t+\delta ,x)|^2|{{\mathcal{F}}_t})]^{\frac{1}{2}}\\
              & \leq&[CE(|X^{t,x;u^{\varepsilon},\beta(u^{\varepsilon})}_{t+\delta} -x|^2|{{\mathcal{F}}_t})]^{\frac{1}{2}},
\end{array}
$$
and since
$E[|X^{t,x;u^{\varepsilon},\beta(u^{\varepsilon})}_{t+\delta}
-x|^2|{{\mathcal{F}}_t}] \leq C(1+|x|^2) \delta $ we deduce that
$|I^1_\delta| \leq C (1+|x|)\delta^{\frac{1}{2}}$. From the
definition of $G^{t,x;
u^{\varepsilon},\beta(u^{\varepsilon})}_{t,t+\delta}[\cdot]$\ (see
(3.14)) we know that the second term $I^2_\delta $ can be written
as£º
$$
\begin{array}{llll}
I^2_\delta  & = & E[W(t+\delta ,x) +\int^{t+\delta}_t
f(s,X^{t,x;u^{\varepsilon},\beta(u^{\varepsilon})}_s,\widetilde{Y}^{t,x;u^{\varepsilon},\beta(u^{\varepsilon})}_s,
\widetilde{Z}^{t,x;u^{\varepsilon},\beta(u^{\varepsilon})}_s,\int_E\widetilde{K}^{t,x;
u^{\varepsilon},\beta(u^{\varepsilon})}_s(e)\\
&&l(X^{t,x;u^{\varepsilon},\beta(u^{\varepsilon})}_s,e)\lambda(de),
             u^{\varepsilon}_s, \beta_s(u^{\varepsilon}_.)) ds   \\
 &         &   -\int^{t+\delta}_t \widetilde{Z}^{t,x; u^{\varepsilon},\beta(u^{\varepsilon})}_s dB_s
 -\int^{t+\delta}_t\int_E\widetilde{K}^{t,x;
u^{\varepsilon},\beta(u^{\varepsilon})}_s(e)\widetilde{\mu}(ds,de)|{{\mathcal{F}}_t}]
  -W(t+\delta ,x)  \\
 &    =  &  E[\int^{t+\delta}_t f(s,X^{t,x; u^{\varepsilon},\beta(u^{\varepsilon})}_s,\widetilde{Y}^{t,x; u^{\varepsilon},
 \beta(u^{\varepsilon})}_s, \widetilde{Z}^{t,x;
 u^{\varepsilon},\beta(u^{\varepsilon})}_s,\int_E\widetilde{K}^{t,x;
u^{\varepsilon},\beta(u^{\varepsilon})}_s(e)l(X^{t,x;u^{\varepsilon},\beta(u^{\varepsilon})}_s,e)\lambda(de),\\
&& u^{\varepsilon}_s,\beta_s(u^{\varepsilon}_.))
 ds|{{\mathcal{F}}_t}],
\end{array}
$$
where
$(\widetilde{Y}^{t,x;u^{\varepsilon},\beta(u^{\varepsilon})}_s,
\widetilde{Z}^{t,x;u^{\varepsilon},\beta(u^{\varepsilon})}_s,
\widetilde{K}^{t,x;
u^{\varepsilon},\beta(u^{\varepsilon})}_s)_{t\leq s\leq t+\delta}$\
is the solution of BSDE (3.15) with the terminal condition
$\eta=W(t+\delta, x).$\ And with the help of the Schwartz
inequality, the estimates (3.3) and (6.4)-(i)  in the Appendix,
we then have that, for some constant $C\in {\mathbb{R}}$\ not
depending on $t$\ and $\delta$,
$$
\begin{array}{lll}
|I^2_\delta | & \leq \delta^{\frac{1}{2}}
     E[\int^{t+\delta}_t |f(s,X^{t,x;u^{\varepsilon},\beta(u^{\varepsilon})}_s,
     \tilde{Y}^{t,x;u^{\varepsilon},\beta(u^{\varepsilon})}_s,\tilde{Z}^{t,x;u^{\varepsilon},\beta(u^{\varepsilon})}_s,
     \int_E\widetilde{K}^{t,x;u^{\varepsilon},\beta(u^{\varepsilon})}_s(e)l(X^{t,x;u^{\varepsilon},
     \beta(u^{\varepsilon})}_s,e)\\
&\ \ \ \ \ \ \lambda(de),u^{\varepsilon}_s,\beta_s(u^{\varepsilon}_.))|^2ds|{{\mathcal{F}}_t}]^{\frac{1}{2}}  \\
& \leq\delta^{\frac{1}{2}}E[\int^{t+\delta}_t
(|f(s,X^{t,x;u^{\varepsilon},\beta(u^{\varepsilon})}_s,0,0,0,u^{\varepsilon}_s,\beta_s(u^{\varepsilon}_.))|+C|
\tilde{Y}^{t,x;u^{\varepsilon}, \beta(u^{\varepsilon})}_s|
+C|\tilde{Z}^{t,x;u^{\varepsilon},\beta(u^{\varepsilon})}_s|\\
&\ \ \ \ \
+C|\int_E\widetilde{K}^{t,x;u^{\varepsilon},\beta(u^{\varepsilon})}_s(e)l(X^{t,x;u^{\varepsilon},
     \beta(u^{\varepsilon})}_s,e)\lambda(de)|)^2ds|{{\mathcal{F}}_t}]^{\frac{1}{2}}\\
& \leq C\delta^{\frac{1}{2}}E[\int^{t+\delta}_t
(1+|X^{t,x;u^{\varepsilon},\beta(u^{\varepsilon})}_s|^2+|
\tilde{Y}^{t,x;u^{\varepsilon}, \beta(u^{\varepsilon})}_s|^2
+|\tilde{Z}^{t,x;u^{\varepsilon},\beta(u^{\varepsilon})}_s|^2\\
&\ \ \ \ \ \ \ \ \ \ \ \
+\int_E|\widetilde{K}^{t,x;u^{\varepsilon},\beta(u^{\varepsilon})}_s(e)|^2\lambda(de))ds|{{\mathcal{F}}_t}]^{\frac{1}{2}}\\
 & \leq C (1+|x|)\delta^{\frac{1}{2}}.
\end{array}
$$
Hence, from (3.19), $$W(t,x)-W(t+\delta ,x) \leq C
(1+|x|)\delta^{\frac{1}{2}} +\varepsilon,$$ and letting $\varepsilon
\downarrow 0$\ we get the second inequality of (3.18). The proof is
complete.\endpf

\section{\large Viscosity solutions of Isaacs' equations with integral-differential operators}

 \hskip1cm In this section we consider the following second order integral-partial differential equations of
 Isaacs' type \be
 \left \{\begin{array}{ll}
 &\!\!\!\!\! \frac{\partial }{\partial t} W(t,x) +  H^{-}(t, x, W, DW, D^2W)=0,
 \hskip 0.5cm   (t,x)\in [0,T)\times {\mathbb{R}}^n ,  \\
 &\!\!\!\!\!  W(T,x) =\Phi (x), \hskip0.5cm   x \in {\mathbb{R}}^n,
 \end{array}\right.
\ee
and
 \be
 \left \{\begin{array}{ll}
 &\!\!\!\!\! \frac{\partial }{\partial t} U(t,x) +  H^{+}(t, x, U, DU, D^2U)=0,
 \hskip 0.5cm   (t,x)\in [0,T)\times {\mathbb{R}}^n ,  \\
 &\!\!\!\!\!  U(T,x) =\Phi (x), \hskip0.5cm   x \in {\mathbb{R}}^n.
 \end{array}\right.
\ee
Their Hamiltonians are given by
$$\begin{array}{lll}
& H^-(t, x, W, DW, D^2W)= \mbox{sup}_{u \in U}\mbox{inf}_{v \in
V}H(t, x, W, DW, D^2W, u, v)
 \end{array}$$
and
$$\begin{array}{lll}& H^+(t, x, U, DU, D^2U)= \mbox{inf}_{v \in
V}\mbox{sup}_{u \in U}H(t, x, U, DU, D^2U, u, v),
 \end{array}$$
$\mbox{respectively, where}$
$$\begin{array}{lll}
& H(t, x, \Psi, D\Psi, D^2\Psi, u,v)=
\frac{1}{2}tr(\sigma\sigma^{T}(t,
x, u, v) D^2\Psi)+ D\Psi.b(t, x, u, v)\\
 &\ \hskip3cm+\int_E(\Psi(t, x+\gamma(t,x,u,v,e))-\Psi(t,x)-D\Psi(t,x).\gamma(t,x,u,v,e))\lambda(de)\\
 &+
 f(t, x, \Psi(t,x), D\Psi(t,x).\sigma(t,x,u,v),\int_E(\Psi(t, x+\gamma(t,x,u,v,e))-\Psi(t,x))l(x,e)\lambda(de), u,
 v),
 \end{array}$$
$\Psi=W \mbox{or}\ U,\ \mbox{resp.,}\ (t, x, u, v)\in [0,
T]\times{\mathbb{R}}^n\times U\times V$. Here the functions $b,
\sigma, f\ \mbox{and}\ \Phi$\ are supposed to satisfy (H3.1) and
(H3.2), respectively.

 In this section we want to prove that the lower value function $W(t, x)$ introduced by (3.9) is the
 viscosity solution of equation (4.1), while the upper value
function $U(t, x)$ defined by (3.10) is the viscosity solution of
equation (4.2). For this we translate Peng's BSDE approach~\cite{P1}
developed in the framework of stochastic control theory driven by
Brownian motion into that of the stochastic differential games
driven by Brownian motion and Poisson random measure. Uniqueness of
the viscosity solution will be shown in the next section for the
class of continuous functions satisfying some growth assumption
which is weaker than the polynomial growth condition. We first
recall the definition of a viscosity solution of equation (4.1). The
definition is analogous for equation (4.2). The reader more
interested in viscosity solutions is referred to Crandall, Ishii and
Lions~\cite{CIL}.

\br We should assume here that there exists a constant C such that
 $$ |\rho(e)|\leq C(1\wedge |e|),\ \mbox{for all}\ e\in E.\eqno{\mbox{(H4.1)}}$$
This assumption (H4.1) is only necessary for the definition 4.1 of
the viscosity solution. If Definition 4.1 is restricted to $W$\ of
linear growth then assumption (H4.1) is not necessary. Notice that
$W$\ defined by (3.9) has this linear growth property. \er

 \bde\mbox{ } A real-valued
continuous function $W\in C([0,T]\times {\mathbb{R}}^n )$ is called \\
  {\rm(i)} a viscosity subsolution of equation (4.1) if $W(T,x) \leq \Phi (x),\ \mbox{for all}\ x \in
  {\mathbb{R}}^n$, and if for all functions $\varphi \in C^3_{l, b}([0,T]\times
  {\mathbb{R}}^n)$ and $(t,x) \in [0,T) \times {\mathbb{R}}^n$ such that $W-\varphi $\ attains
  a local maximum at $(t, x)$,
    \be
    \begin{array}{lll}
     &\frac{\partial \varphi}{\partial t} (t,x) +\mbox{sup}_{u \in U}\mbox{inf}_{v \in
V}\{A^{u,v}\varphi(t,x)+ B^{\delta,u,v}(W,\varphi)(t,x)\\
 &\ \ \ \ \ +f(t, x,
W(t,x),D\varphi(t,x).\sigma(t,x,u,v),C^{\delta,u,v}(W,\varphi)(t,x),
u, v)\}\geq 0, \end{array}\ee for any $\delta>0,$\ where
$$A^{u,v}\varphi(t,x)=\frac{1}{2}tr(\sigma\sigma^{T}(t, x,
 u, v) D^2\varphi(t,x))+ D\varphi(t,x).b(t, x, u, v),$$
$$
\begin{array}{lll}
B^{\delta,u,v}(W,\varphi)(t,x)&=\int_{E_\delta}(\varphi(t,
x+\gamma(t,x,u,v,e))-\varphi(t,x)-D\varphi(t,x).\gamma(t,x,u,v,e))\lambda(de)\\
&+\int_{E_{\delta}^c}(W(t,
x+\gamma(t,x,u,v,e))-W(t,x)-D\varphi(t,x).\gamma(t,x,u,v,e))\lambda(de),
\end{array}
$$
and
$$
\begin{array}{lll}
C^{\delta,u,v}(W,\varphi)(t,x)&=&\int_{E_\delta}(\varphi(t,
x+\gamma(t,x,u,v,e))-\varphi(t,x))l(x,e)\lambda(de)\\
&& +\int_{E_{\delta}^c}(W(t,
x+\gamma(t,x,u,v,e))-W(t,x))l(x,e)\lambda(de),
\end{array}
$$
with $E_\delta=\{e\in E||e|<\delta\}$.

\noindent{\rm(ii)} a viscosity supersolution of equation (4.1) if
$W(T,x) \geq \Phi (x), \mbox{for all}\ x \in
  {\mathbb{R}}^n$, and if for all functions $\varphi \in C^3_{l, b}([0,T]\times
  {\mathbb{R}}^n)$ and $(t,x) \in [0,T) \times {\mathbb{R}}^n$\ such that $W-\varphi $\ attains
  a local minimum at $(t, x)$,
      \be
    \begin{array}{lll}
     &\frac{\partial \varphi}{\partial t} (t,x) +\mbox{sup}_{u \in U}\mbox{inf}_{v \in
V}\{A^{u,v}\varphi(t,x)+ B^{\delta,u,v}(W,\varphi)(t,x)\\
 &\ \ \ \ \ +f(t, x,
W(t,x),D\varphi(t,x).\sigma(t,x,u,v),C^{\delta,u,v}(W,\varphi)(t,x),
u, v)\}\leq 0. \end{array}\ee
 {\rm(iii)} a viscosity solution of equation (4.1) if it is both a viscosity sub- and a supersolution of equation
     (4.1).\ede
\br \mbox{  }$C^3_{l, b}([0,T]\times {\mathbb{R}}^n)$ denotes the
set of  real-valued functions that are continuously
differentiable up to the third order and whose derivatives of
order from 1 to 3 are bounded.\er

\noindent In analogy to ~\cite{BBP} we have the following result:

\bl In the definition of W being a viscosity sub- (resp.,
super-)solution of (4.1), we can replace
$$\begin{array}{rcl}
B^{\delta,u,v}(W,\varphi)(t,x)&=&B^{u,v}\varphi(t,x),\\
C^{\delta,u,v}(W,\varphi)(t,x)&=&C^{u,v}\varphi(t,x),
\end{array}
$$
where
$$\begin{array}{lll}
&B^{u,v}\varphi(t,x)=\int_{E}(\varphi(t,
x+\gamma(t,x,u,v,e))-\varphi(t,x)-D\varphi(t,x).\gamma(t,x,u,v,e))\lambda(de),\\
&C^{u,v}\varphi(t,x)=\int_{E}(\varphi(t,
x+\gamma(t,x,u,v,e))-\varphi(t,x))l(x,e)\lambda(de).
\end{array}
$$
\el \noindent Proof: We only consider the subsolution case, the supersolution case can be treated analogously.\\
If $(t,x) \in [0,T) \times {\mathbb{R}}^n$\ such that $W-\varphi $\
attains a global maximum at $(t, x)$\ we have
$W(s,y)-\varphi(s,y)\leq W(t,x)-\varphi(t,x),\ \mbox{for all}\ (s,y)
  \in[0,T]\times {\mathbb{R}}^n$. Therefore, $W(t,y)-W(t,x)\leq \varphi(t,y) -\varphi(t,x),\ \ \mbox{for
  any}\ \ y\in{\mathbb{R}}^n$\ and this yields, for any $\delta>0,$
  $$\begin{array}{rcl}
B^{\delta,u,v}(W,\varphi)(t,x)&\leq &B^{u,v}\varphi(t,x),\\
C^{\delta,u,v}(W,\varphi)(t,x)&\leq &C^{u,v}\varphi(t,x).
\end{array}
$$
Because $f$ is increasing in $k$, from (4.3) we get  \be
    \begin{array}{lll}
     &\frac{\partial \varphi}{\partial t} (t,x) +\mbox{sup}_{u \in U}\mbox{inf}_{v \in
V}\{A^{u,v}\varphi(t,x)+ B^{u,v}\varphi(t,x)\\
 &\ \ \ \ \ +f(t, x,
W(t,x), D\varphi.\sigma(t,x,u,v),C^{u,v}\varphi(t,x), u, v)\}\geq
0.\end{array}\ee It remains to show this last condition (4.5)
implies (4.3). Changing $\varphi$\ into
$\varphi-(\varphi(t,x)-W(t,x))$, we may assume that
$W(t,x)=\varphi(t,x)$. Then $W(s,y)\leq \varphi(s,y),\ \mbox{for
all}\ (s,y)\in[0,T]\times {\mathbb{R}}^n$. Moreover, we may assume
without loss of generality that,\\
 $$\begin{array}{ll}
 &\mbox{(i)}\mbox{ for all}\ \alpha>0, \mbox{there exists some}\
 \eta_{\alpha}>0,\ \mbox{with}\ \eta_{\alpha}\rightarrow0\ \mbox{as}\
 \alpha\rightarrow0,\ \mbox{such that, for all}\ \ \ \ \
\ \\
 &\ \ \ \ (s,y)\in[0,T]\times {\mathbb{R}}^n\ \mbox{with}\
|(s,y)-(t,x)|>\alpha,\ \varphi(s,y)-W(s,y)\geq \eta_{\alpha}.\ \ \ \
\
\end{array}
 $$
Furthermore, there exists a sequence of elements $\varphi_\alpha$ in
$C_{l,b}^3([0,T]\times {\mathbb{R}}^n)$\ with the following
properties:
  $$
  \begin{array}{ll}
  &\mbox{(ii)}\ \varphi_\alpha(s,y)=\varphi(s,y)\geq W(s,y),\mbox{ if  }
                                 |(s,y)-(t,x)|\notin(\alpha,\frac{1}{\alpha});\\[.2cm]
  &\mbox{(iii)}\ \varphi_\alpha(s,y)\geq W(s,y),\mbox{ if }
                                \ \alpha\leq|(s,y)-(t,x)|\leq\frac{1}{\alpha};\\[.2cm]
  &\mbox{(iv)}\ \varphi_\alpha(s,y)\leq W(s,y)+ \eta_{\alpha}, \mbox{ if }3\alpha\leq
                        |(s,y)-(t,x)|\leq \frac{1}{\alpha}-2\alpha;\\[.2cm]
&\mbox{(v)}\ \varphi_\alpha(s,y)\leq \varphi(s,y),\ \mbox{for all }\
(s,y)\in [0,T]\times {\mathbb{R}}^n;\\[.2cm]
&\mbox{(vi)}\ \mbox{There exists some }\ \bar{\rho}_\alpha>0\
\mbox{with}\ \bar{\rho}_\alpha\rightarrow 0\ (\alpha\downarrow 0)\
\mbox{such that
}\\
&\ \ \ \ 0\leq\varphi_\alpha(s,y)-W(s,y)\leq \bar{\rho}_\alpha,\
\mbox{for all}\ (s, y)\ \mbox{satisfying}\ |(s,y)-(t,x)|\leq \frac{1}{\alpha}-2\alpha.\\[.2cm]
 \end{array}
 $$
Then, obviously, we have
  $D\varphi_\alpha(t,x)=D\varphi(t,x),
  \frac{\partial\varphi_\alpha(t,x)}{\partial t}=
 \frac{\partial\varphi(t,x)}{\partial t}, D^2\varphi_\alpha(t,x)=D^2\varphi(t,x).$\ Thus, since
  $\varphi_\alpha(t,x)=W(t,x)$\ and $\varphi_\alpha(s,y)\geq
  W(s,y)$, it follows from (4.5) that
   $$
     \begin{array}{ll}
      &\frac{\partial \varphi_\alpha(t,x)}{\partial t}+\mbox{sup}_{u \in U}\mbox{inf}_{v \in
V}\{A^{u,v}\varphi_\alpha(t,x)+ B^{u,v}\varphi_\alpha(t,x)\\
 &\ \ \ \ \ +f(t, x,
W(t,x), D\varphi_\alpha.\sigma(t,x,u,v), C^{u,v}\varphi_\alpha(t,x),
u, v)\}\geq 0.
     \end{array}
     $$
From the above property (v) and the monotonicity of $f$, we get
 $$
     \begin{array}{ll}
      &\frac{\partial \varphi(t,x)}{\partial t}+\mbox{sup}_{u \in U}\mbox{inf}_{v \in
V}\{A^{u,v}\varphi(t,x)+ B^{\delta,u,v}(\varphi_\alpha, \varphi)(t,x)\\
 &\ \ \ \ \ +f(t, x,
W(t,x), D\varphi.\sigma(t,x,u,v), C^{\delta, u,v}(\varphi_\alpha,
\varphi)(t,x), u, v)\}\geq 0.
     \end{array}
     $$
Finally, by (vi) and the Lebesgue dominated convergence theorem we
deduce that
$$\lim\limits_{\alpha\rightarrow0}
  \mbox{sup}_{u \in U}\mbox{inf}_{v \in
V} B^{\delta,u,v}(\varphi_\alpha, \varphi)(t,x)=\mbox{sup}_{u \in
U}\mbox{inf}_{v \in V} B^{\delta,u,v}(W, \varphi)(t,x);$$
 $$\lim\limits_{\alpha\rightarrow0}
  \mbox{sup}_{u \in U}\mbox{inf}_{v \in
V} C^{\delta,u,v}(\varphi_\alpha, \varphi)(t,x)=\mbox{sup}_{u \in
U}\mbox{inf}_{v \in V} C^{\delta,u,v}(W, \varphi)(t,x).$$

\noindent Indeed, since $\varphi_\alpha(t, .),\ W(t, .)$\ are
continuous and coincide in $x$\ we get
 $$
 \begin{array}{ll}
&|B^{\delta,u,v}(\varphi_\alpha,
\varphi)(t,x)-B^{\delta,u,v}(W,\varphi)(t,x)|\\
\leq &\int_{E_{\delta}^c}|\varphi_\alpha(t,
x+\gamma(t,x,u,v,e))-W(t,x+\gamma(t,x,u,v,e))|\lambda(de)\\
\leq & \bar{\rho}_{\alpha}\lambda(E_{\delta}^c)+C\int_{E_{\delta}^c}
\textbf{I}_{\{|\gamma(t,x,u,v,e)|>\frac{1}{\alpha}-2\alpha\}}\lambda(de)\\
& \longrightarrow0\ (\alpha\rightarrow 0)\ \mbox{uniformly in}\
(u,v)\in U\times V.
     \end{array}
     $$
\noindent The second convergence uses the same argument. Therefore,
letting $\varepsilon\rightarrow0$ in the above estimate
  yields the desired result.\endpf

\noindent\br\mbox{    }As concerns the construction of the sequence $(\varphi_\alpha)$, we also refer the reader to Remark 4.3 in Li and Peng~\cite{LP}.\er\vskip0.1cm

Let us now first prove that the lower value function $W(t,x)$ is a
viscosity solution of equation (4.1).
 \bt Under the assumptions (H3.1) and (H3.2) the lower value function $W(t,x)$ is a viscosity solution of
equation (4.1). \et

\noindent For the proof of this theorem we need four auxiliary
lemmas. To abbreviate notation we put, for some arbitrarily
chosen but fixed $\varphi \in C^3_{l, b} ([0,T] \times
{\mathbb{R}}^n)$,
 \be
\begin{array}{lll}
     &\ \hskip2cm F(s,x,y,z, k, u, v)= \frac{\partial }{\partial s}\varphi (s,x) +
     A^{u,v}\varphi(s,x)+B^{u,v}\varphi(s,x)\\
        &+ f(s, x, y+\varphi (s,x), z+ D\varphi (s,x).\sigma(s,x,u, v), \int_Ek(e)l(x,e)\lambda(de)+C^{u,v}\varphi(s,x),u, v), \\
     \end{array}
\ee $(s,x,y,z,k,u, v)\in [0,T] \times {\mathbb{R}}^n \times
{\mathbb{R}} \times {\mathbb{R}}^d\times L^{2}(E, {\cal{B}}(E),
\lambda;{\mathbb{R}}) \times U \times V,$\ and we consider the
following BSDE defined on the interval $[t,t+\delta]\ (0<\delta\leq
T-t):$
    \be \left \{\begin{array}{rl}
      -dY^{1,u,v}_s =&\!\!\!\! F(s,X^{t,x;u,v}_s, Y^{1,u,v}_s, Z^{1,u,v}_s,K^{1,u,v}_s, u_s,v_s)ds
                   -Z^{1,u,v}_s dB_s-\int_EK^{1,u,v}_s(e)\widetilde{\mu}(ds,de), \\
     Y^{1,u,v}_{t+\delta}=&\!\!\!\!0,
     \end{array}\right.
     \ee
     where the process $X^{t,x,u,v}$\ has been introduced by equation
     $(3.1)$\ and $u(\cdot) \in {\mathcal{U}}_{t, t+\delta},\ v(\cdot) \in
     {\mathcal{V}}_{t, t+\delta}$.
\br\mbox{}It is not hard to check that $F(s,X^{t,x;u,v}_s,
y,z,k,u_s,v_s)$\ satisfies (A1) and (A2). Thus, due to Lemma 2.1,
equation (4.7) has a unique solution. \er

      We can characterize the solution process $Y^{1,u,v}$ as follows:
    \bl\mbox{  } For every $s\in [t,t+\delta]$, we have the following
relationship:
    \be
     Y^{1,u,v}_s = G^{t,x;u,v}_{s,t+\delta} [\varphi (t+\delta ,X^{t,x;u,v}_{t+\delta})]
                -\varphi (s,X^{t,x;u,v}_s), \hskip 0.5cm
              \mbox{ P-a.s.} \ee
 \el
\noindent \textbf{Proof}: We recall that $G^{t,x;u,v}_{s,t+\delta}
[\varphi (t+\delta, X^{t,x;u,v}_{t+\delta})]$ is defined with the
help of the solution of the BSDE
     $$
     \left \{\begin{array}{rl}
     -dY^{u,v}_s =\!\!\! & f(s, X^{t,x;u,v}_s, Y^{u,v}_s, Z^{u,v}_s,
     \int_EK^{u,v}_s(e)l(X^{t,x;u,v}_s,e)\lambda(de), u_s,v_s)ds
      -Z^{u,v}_s dB_s\\
      &-\int_EK^{u,v}_s(e)\widetilde{\mu}(ds,de),  \hskip 3cm s\in [t,t+\delta], \\
     Y^{u,v}_{t+\delta}=\!\!\!& \varphi (t+\delta
     ,X^{t,x;u,v}_{t+\delta}),
     \end{array}\right.
     $$
by the following formula:
     \be
     G^{t,x;u,v}_{s,t+\delta} [\varphi (t+\delta ,X^{t,x;u,v}_{t+\delta})] =Y^{u,v}_s, \hskip 0.5cm
     s\in [t,t+\delta]  \ee
(see (3.14)). Therefore, we only need to prove that
$Y^{u,v}_s-\varphi (s,X^{t,x;u,v}_s)\equiv Y^{1,u,v}_s.$\ This
result can be obtained easily by applying It$\hat{o}$'s formula to
$\varphi (s,X^{t,x;u,v}_s)$. Indeed, we get that the stochastic
differentials of $Y^{u,v}_s -\varphi (s,X^{t,x;u,v}_s)$\ and
$Y^{1,u,v}_s$\ coincide, while at the terminal time $t+\delta$,
$Y^{u,v}_{t+\delta} - \varphi (t+\delta ,X^{t,x;u,v}_{t+\delta}) =0
= Y^{1,u,v}_{t+\delta}.$
     So the proof is complete.\endpf \vskip 0.3cm

Now we consider the following simple BSDE in which the driving
process $X^{t,x;u,v}$ is replaced by its deterministic initial
value $x$: \be
     \left \{\begin{array}{rl}
     -dY^{2,u,v}_s=&\!\!\!  F(s,x,Y^{2,u,v}_s ,Z^{2,u,v}_s, K^{2,u,v}_s, u_s, v_s)ds -Z^{2,u,v}_s dB_s-
     \int_EK^{2,u,v}_s(e)\widetilde{\mu}(ds,de),
         \\
     Y^{2,u,v}_{t+\delta}=&\!\!\! 0,
         \hskip  3cm s\in [t,t+\delta],
     \end{array}\right.
 \ee
where $u(\cdot) \in {\mathcal{U}}_{t, t+\delta},\ v(\cdot) \in
     {\mathcal{V}}_{t, t+\delta}$.
The following lemma will allow us to neglect the difference
$|Y^{1,u,v}_t-Y^{2,u,v}_t|$ for sufficiently small $\delta >0$.

\bl For every $u \in {\mathcal{U}}_{t, t+\delta},\ v \in
{\mathcal{V}}_{t, t+\delta},$\ we have
 \be |Y^{1,u,v}_t-Y^{2,u,v}_t| \leq C\delta^{\frac{3}{2}},\ \ \mbox{P-a.s.},
       \ee
 where C is independent
 of the control processes $u$\ and $v$.
\el
 \noindent \textbf{Proof}: From Proposition 1.1 in~\cite{BBP}, we have for all $p\geq 2$\ the existence
 of some $C_{p}\in {\mathbb{R}}^+$\ such that
      \be
       E [\sup \limits_{t\leq s \leq t+\delta} |X^{t,x;u,v}_s -x|^p|{{\mathcal{F}}_t}] \leq
       C_p\delta(1+|x|^p),\ \
    \mbox{P-a.s., \ uniformly in}\ u \in {\mathcal{U}}_{t, t+\delta}, v \in
{\mathcal{V}}_{t, t+\delta}.
      \ee
   We now apply Lemma 2.3 combined with (4.12) to equations (4.7) and (4.10). For this we
set in Lemma 2.3:
       $$\xi_1 =\xi_2 =0,\ g(s,y,z) =F(s,X^{t,x,u,v}_s,y,z,k,u_s,v_s),$$
       $$\varphi_1(s)=0,\ \varphi_2(s)=F(s,x,Y_s^{2, u, v},Z_s^{2, u, v},K_s^{2, u, v},u_s,v_s)
       -F(s,X^{t,x,u,v}_s,Y_s^{2, u, v},Z_s^{2, u, v},K_s^{2, u, v},u_s,v_s).
       $$
       Obviously, the function $g$\ is Lipschitz with respect to
       $(y,z,k)$. We also notice that
       $$\begin{array}{lll}
&B^{u,v}\varphi(s,x)=\int_{E}(\varphi(s,
x+\gamma(s,x,u,v,e))-\varphi(s,x)-D\varphi.\gamma(s,x,u,v,e))\lambda(de)\\
&\ \hskip2cm
=\int_E\int_0^1(1-\theta)tr(D^2\varphi(s,x+\theta\gamma(s,x,u,v,e))\gamma\gamma^T(s,x,u,v,e))d\theta\lambda(de);\\
&C^{u,v}\varphi(s,x)=\int_{E}(\varphi(s,
x+\gamma(s,x,u,v,e))-\varphi(s,x))l(x,e)\lambda(de)\\
&\ \hskip2cm
=\int_E\int_0^1D\varphi(s,x+\theta\gamma(s,x,u,v,e))\gamma(s,x,u,v,e)d\theta
l(x,e)\lambda(de).
\end{array}
$$
 Then, we can get $|\varphi_2(s)|\leq C(1+|x|^2)(|X^{t,x;u,v}_s -x|+|X^{t,x;u,v}_s -x|^3),$\
       for $s\in [t, t+\delta], (t, x)\in [0, T)\times {\mathbb{R}}^n$, $u\in {\mathcal{U}}_{t, t+\delta}, v \in {\mathcal{V}}_{t, t+\delta}.
       $
       Thus, with the notation $\rho_0 (r) =(1+|x|^2)(r+r^3),\ r\geq 0, $\ we have
      $$
       \begin{array}{lll}
      & E[\int^{t+\delta}_t (|Y^{1,u,v}_s -Y^{2,u,v}_s|^2 +|Z^{1,u,v}_s -Z^{2,u,v}_s|^2)ds|{\mathcal{F}}_{t}]  \\
       &\ \hskip2cm + E[\int^{t+\delta}_t\int_E |K^{1,u,v}_s(e) -K^{2,u,v}_s(e)|^2
       \lambda(de)ds|{\mathcal{F}}_{t}]\\
       & \leq  CE[\int^{t+\delta}_t \rho^2_0 (|X^{t,x,u,v}_s -x|) ds|{\mathcal{F}}_{t}]\\
       & \leq C \delta E[\sup \limits_{t\leq s \leq t+\delta}\rho^2_0 (|X^{t,x,u,v}_s
       -x|)|{\mathcal{F}}_{t}]\\
       &\leq C\delta^2.
       \end{array}
     $$
Therefore,
       $$
       \begin{array}{llll}
      && |Y^{1,u,v}_t -Y^{2,u,v}_t|  =|E[(Y^{1,u,v}_t -Y^{2,u,v}_t )|{\mathcal{F}}_{t}]|  \\
       & = & |E[\int^{t+\delta}_t (F(s,X^{t,x,u,v}_s,Y^{1,u,v}_s,Z^{1,u,v}_s,K^{1,u,v}_s,u_s,
       v_s)\\
        &&\ \hskip2cm         -F(s,x,Y^{2,u,v}_s ,Z^{2,u,v}_s,K^{2,u,v}_s,u_s, v_s)) ds|{\mathcal{F}}_{t}]|   \\
       & \leq & CE [\int^{t+\delta}_t (\rho_0 (|X^{t,x,u,v}_s -x|) +|Y^{1,u,v}_s -Y^{2,u,v}_s|
       +|Z^{1,u,v}_s -Z^{2,u,v}_s|)ds|{\mathcal{F}}_{t}]\\
      &&\ \hskip1cm +CE [\int^{t+\delta}_t|\int_E(K^{1,u,v}_s(e)-K^{2,u,v}_s(e))l(x,e)\lambda(de)|ds|{\mathcal{F}}_{t}]  \\
       & \leq & CE[\int^{t+\delta}_t\rho_0 (|X^{t,x,u,v}_s -x|)ds|{\mathcal{F}}_{t}] + C\delta^{\frac{1}{2}} \{ E[\int^{t+\delta}_t
            |Y^{1,u,v}_s-Y^{2,u,v}_s|^2ds|{\mathcal{F}}_{t}]^{\frac{1}{2}}\\
         &&   + E[\int^{t+\delta}_t|Z^{1,u,v}_s -Z^{2,u,v}_s|^2ds|{\mathcal{F}}_{t}]^{\frac{1}{2}}
         + E[\int^{t+\delta}_t\int_E|K^{1,u,v}_s(e) -K^{2,u,v}_s(e)|^2\lambda(de)ds|{\mathcal{F}}_{t}]^{\frac{1}{2}}\}  \\
       & \leq & C\delta^{\frac{3}{2}}.
       \end{array}
       $$
   Thus, the proof is complete.\endpf
     \vskip 0.3cm
\bl Let $Y_0 (\cdot)$ be the solution of the following ordinary
differential equation:
     \be
     \left \{\begin{array}{lll}
     - {\dot{Y}}_0 (s)  & =  & F_0(s,x,Y_0 (s),0,0),\ \ s\in [t,t+\delta], \\
      Y_0 (t+\delta )   & =  & 0,
     \end{array}\right.
    \ee
  where the function $F_0$ is defined by
    \be\begin{array}{lll}
     &\ \hskip1cm F_0(s,x,y,z,k)=\mbox{sup}_{u\in U}\mbox{inf}_{v \in
V} F(s,x,y,z,k,u,v),\\
&(s,x,y,z,k)\in [0, T]\times {\mathbb{R}}^n \times
{\mathbb{R}}\times {\mathbb{R}}^d
     \times L^{2}(E, {\cal{B}}(E),
\lambda;{\mathbb{R}}).
\end{array}
 \ee Then, P-a.s.,\be\mbox{esssup}_{u \in {\mathcal{U}}_{t,
t+\delta}}\mbox{essinf}_{v \in {\mathcal{V}}_{t, t+\delta}}
Y^{2,u,v}_t =Y_0 (t). \ee
  \el
 \noindent \textbf{Proof}: Obviously,
 $F_0(s,x,y,z,k)$\ is Lipschitz in $(y, z,k)$, uniformly with respect to $(s, x).$\ This guarantees the existence and uniqueness for equation
 (4.13). We first introduce the function \be\begin{array}{lll}
     &F_1(s,x,y,z,k,u)=\mbox{inf}_{v\in V} F(s,x,y,z,k,u,v),\\
    &\ \hskip1cm (s,x,y,z,k,u)\in [0, T]\times {\mathbb{R}}^n \times {\mathbb{R}}\times {\mathbb{R}}^d
     \times L^{2}(E, {\cal{B}}(E),
\lambda;{\mathbb{R}})\times U, \end{array}\ee and consider the BSDE
\be
     \left \{\begin{array}{lll}
     - d{Y}^{3,u}_s  & =  & F_1(s,x,{Y}^{3,u}_s,{Z}^{3,u}_s,{K}^{3,u}_s,u_s)ds-{Z}^{3,u}_sdB_s
     -\int_E{K}^{3,u}_s(e)\widetilde{\mu}(ds,de), \\
      {Y}^{3,u}_{t+\delta}   & =  & 0,\ \ \ s\in [t,t+\delta], \ \
       \end{array}\right.
    \ee
$\mbox{for}\ u\in {\mathcal{U}}_{t, t+\delta}$. We notice that since
$F_1(s,x,y,z,k,u_s)$\ is Lipschitz in $(y, z,k)$, for every $u\in
{\mathcal{U}}_{t, t+\delta},$\ there exists a unique solution
$({Y}^{3,u}, {Z}^{3,u}, {K}^{3,u})$\ to the BSDE (4.17). Moreover,
$${Y}^{3,u}_t= \mbox{essinf}_{v(\cdot) \in {\mathcal{V}}_{t, t+\delta}}
Y^{2,u,v}_t,\ \mbox{P-a.s.},\ \mbox{for all}\ u\in
      {\mathcal{U}}_{t, t+\delta}.$$
Indeed, from the definition of $F_1$\ and Lemma 2.2 (comparison
theorem) we have
$${Y}^{3,u}_t\leq \mbox{essinf}_{v(\cdot) \in {\mathcal{V}}_{t, t+\delta}}
Y^{2,u,v}_t,\ \mbox{P-a.s.},\ \mbox{for all}\ u\in
      {\mathcal{U}}_{t, t+\delta}.$$
On the other hand, there exists a measurable function $v^3: [t,
T]\times
{\mathbb{R}}^n\times{\mathbb{R}}\times{\mathbb{R}}^d\times{\mathbb{R}}\times
U\rightarrow V$\ such that
$$F_1(s,x,y,z,k,u)= F(s,x,y,z,k,u,v^3(s,x,y,z,k,u)),\ \mbox{for any}\ s, x, y, z,k, u.$$
Then, given an arbitrary $u\in {\cal{U}}_{t, t+\delta}$\ we put
$$\widetilde{v}_s^3:=
     v^{3} (s, x, {Y}^{3,u}_s, {Z}^{3,u}_s, {K}^{3,u}_s, u_s ), \  s\in [t,
     t+\delta],
$$
and observe that $\widetilde{v}^3 \in {\mathcal{V}}_{t, t+\delta},$\
and
$$F_1(s,x,{Y}^{3,u}_s, {Z}^{3,u}_s,{K}^{3,u}_s, u_s)= F(s,x,{Y}^{3,u}_s, {Z}^{3,u}_s, {K}^{3,u}_s,u _s, \widetilde{v}_s^3),\ s\in [t, t+\delta].$$
Consequently, from the uniqueness of the solution of the BSDE it
follows that $$({Y}^{3,u},
{Z}^{3,u},{K}^{3,u})=({Y}^{2,u,\widetilde{v}^3}, {Z}^{2,u,
\widetilde{v}^3}, {K}^{2,u, \widetilde{v}^3})$$ and, in particular,
${Y}^{3,u}_t={Y}^{2,u,\widetilde{v}^3}_t,\ \mbox{P-a.s.}$\ This
proves that
$${Y}^{3,u}_t= \mbox{essinf}_{v\in {\mathcal{V}}_{t, t+\delta}}
Y^{2,u,v}_t,\ \mbox{P-a.s.},\ \mbox{for all}\ u\in
      {\mathcal{U}}_{t, t+\delta}.$$
 Finally, since $F_0(s,x,y,z,k)=\mbox{sup}_{u \in
U}F_1(s,x,y,z,k,u),$ \ an argument similar to that developed above
yields
$$Y_0 (t)=\mbox{esssup}_{u \in
{\mathcal{U}}_{t, t+\delta}} Y^{3,u}_t(=\mbox{esssup}_{u \in
{\mathcal{U}}_{t, t+\delta}}\mbox{essinf}_{v \in {\mathcal{V}}_{t,
t+\delta}} Y^{2,u,v}_t), \ \mbox{P-a.s.}$$ It uses the fact that
equation (4.13) can be considered as a BSDE with solution $(Y_s,
Z_s, K_s)= (Y_0(s), 0,0).$
 The proof is complete.\endpf
\bl For every $u \in {\mathcal{U}}_{t, t+\delta}, v \in
{\mathcal{V}}_{t, t+\delta},$\ we have\be\begin{array}{ll}
&E[\int_t^{t+\delta}|Y^{2,u,v}_s|ds|{\mathcal{F}}_{t}]+E[\int_t^{t+\delta}|Z^{2,u,v}_s|ds|{\mathcal{F}}_{t}]
+E[\int_t^{t+\delta}|\int_EK^{2,u,v}_s(e)l(x,e)\lambda(de)|ds|{\mathcal{F}}_{t}]\\
&\leq C\delta^{\frac{3}{2}},\ \mbox{P-a.s.},\end{array} \ee  where
the constant $C$\ is independent of\ $t,\ \delta$\ and the control
processes $u,\ v$.\el
 \noindent \textbf{Proof}: Since $F(s, x, \cdot, \cdot,\cdot, u, v)$\ has a linear
 growth in $(y, z,k)$, uniformly in $(u, v)$, we get from Lemma 2.3 that, for some constant $C$\ independent of $\delta$\
 and the control processes $u,\ v,\ P\mbox{-a.s.},$
 $$|Y^{2,u,v}_s|^2\leq C\delta,\ E[\int_s^{t+\delta}|Z^{2,u,v}_r|^2dr|{\cal{F}}_s]\leq
C\delta,$$
$$E[\int_s^{t+\delta}\int_E|K^{2,u,v}_r(e)|^2\lambda(de)dr|{\cal{F}}_s]\leq
C\delta,\ t\leq s\leq t+\delta.$$ On the other hand, from equation
(4.10),
 $$\begin{array}{ll}
 |Y^{2,u,v}_s|&\leq E[\int_s^{t+\delta}|F(r,x,{Y}^{2, u, v}_r, {Z}^{2, u, v}_r,{K}^{2, u, v}_r, u _r, v_r)|dr|{\cal{F}}_s]\\
 &\leq CE[\int_s^{t+\delta}(1+ |x|^2+|{Y}^{2, u, v}_r|+ |{Z}^{2, u, v}_r|+|\int_EK^{2,u,v}_r(e)l(x,e)\lambda(de)|)dr|{\cal{F}}_s]\\
 &\leq C\delta+C\sqrt{\delta}(E[\int_s^{t+\delta}|Z^{2,u,v}_r|^2dr|{\cal{F}}_s])^{\frac{1}{2}}
 +C\sqrt{\delta}(E[\int_s^{t+\delta}\int_E|K^{2,u,v}_r(e)|^2\lambda(de)dr|{\cal{F}}_s])^{\frac{1}{2}}\\
 &\leq C\delta,\ \mbox{P-a.s.},\ \ s\in[t,
 t+\delta],
 \end{array}$$
and, applying It\^o formula to $|Y^{2,u,v}_s|^2$ we can get
$$E[\int_t^{t+\delta}|Z^{2,u,v}_s|^2ds|{\cal{F}}_t]+E[\int_t^{t+\delta}\int_E|K^{2,u,v}_s(e)|^2\lambda(de)ds|{\cal{F}}_t]\leq
C\delta^2, \ \ P\mbox{-a.s.}$$\ Finally,
$$\begin{array}{ll}
&E[\int_t^{t+\delta}|Y^{2,u,v}_s|ds|{\cal{F}}_t]+E[\int_t^{t+\delta}|Z^{2,u,v}_s|ds|{\cal{F}}_t]
+E[\int_t^{t+\delta}|\int_EK^{2,u,v}_s(e)l(x,e)\lambda(de)|ds|{\mathcal{F}}_{t}]\\
&\leq
C\delta^2+\delta^{\frac{1}{2}}\{E[\int_t^{t+\delta}|Z^{2,u,v}_s|^2ds|{\cal{F}}_t]\}^{\frac{1}{2}}
+C\delta^{\frac{1}{2}}\{E[\int_t^{t+\delta}\int_E|K^{2,u,v}_s(e)|^2\lambda(de)ds|{\cal{F}}_t]\}^{\frac{1}{2}}\\
&\leq C\delta^{\frac{3}{2}}, \ \ P\mbox{-a.s.}\end{array}$$ The
proof is complete.\endpf \vskip0.2cm
 Now we are able to give the proof of Theorem 4.1:

\noindent \textbf{Proof}: (1) Obviously, $W(T,x)=\Phi (x),\ x\in
{\mathbb{R}}^n $. Let us show in a first step that $W$ is a
viscosity supersolution. For this we suppose that $\varphi \in
C^3_{l,b} ([0,T] \times {\mathbb{R}}^n)$,\ and $(t,x)\in [0,
T)\times {\mathbb{R}}^n$\ are such that $W-\varphi$\ attains a
minimum at $(t,x).$\ Notice that we can replace the condition of a
local minimum by that of a global one in the definition of the
viscosity supersolution since $W$ is continuous and of at most
linear growth. Without loss of generality we may also suppose that
$\varphi (t,x)=W(t,x)$.\ Then, due to the DPP (see Theorem 3.1),
     $$
     \varphi (t,x) =W(t,x) =\mbox{essinf}_{\beta \in {\mathcal{B}}_{t, t+\delta}}\mbox{esssup}_{u \in
{\mathcal{U}}_{t, t+\delta}}G^{t,x;u,\beta(u)}_{t,t+\delta}
[W(t+\delta, X^{t,x;u,\beta(u)}_{t+\delta})],\ 0\leq\delta\leq
T-t,
     $$
 and from $W\geq \varphi$\ and the monotonicity property of $G^{t,x;u,\beta(u)}_{t,t+\delta}[\cdot]$\ (see Lemma 2.2)  we obtain
     $$
     \mbox{essinf}_{\beta \in {\mathcal{B}}_{t, t+\delta}}\mbox{esssup}_{u \in
{\mathcal{U}}_{t, t+\delta}} \{G^{t,x;u,\beta(u)}_{t,t+\delta}
[\varphi(t+\delta, X^{t,x;u,\beta(u)}_{t+\delta})] -\varphi
(t,x)\}\leq 0,\ \mbox{P-a.s.}
     $$
 Thus, from Lemma 4.2,
    $$
      \mbox{essinf}_{\beta \in {\mathcal{B}}_{t, t+\delta}}\mbox{esssup}_{u \in
{\mathcal{U}}_{t, t+\delta}} Y^{1,u,\beta(u)}_t \leq 0,\
\mbox{P-a.s.},
    $$
and further, from Lemma 4.3 we have
     $$
     \mbox{essinf}_{\beta \in {\mathcal{B}}_{t, t+\delta}}\mbox{esssup}_{u \in
{\mathcal{U}}_{t, t+\delta}} Y^{2,u,\beta(u)}_t \leq
C\delta^{\frac{3}{2}},\ \mbox{P-a.s.}
     $$
     Consequently, since $\mbox{essinf}_{v \in
{\mathcal{V}}_{t, t+\delta}} Y^{2,u,v}_t\leq Y^{2,u,\beta(u)}_t, \
\beta\in{\mathcal{B}}_{t, t+\delta} $, we get
$$
\mbox{esssup}_{u \in {\mathcal{U}}_{t, t+\delta}}\mbox{essinf}_{v
\in {\mathcal{V}}_{t, t+\delta}}
Y^{2,u,v}_t\leq\mbox{essinf}_{\beta \in {\mathcal{B}}_{t,
t+\delta}}\mbox{esssup}_{u \in {\mathcal{U}}_{t, t+\delta}}
Y^{2,u,\beta(u)}_t \leq C\delta^{\frac{3}{2}},\ \mbox{P-a.s.},
$$
and Lemma 4.4 implies
     $$
     Y_0 (t) \leq C\delta^{\frac{3}{2}},\ \mbox{P-a.s.},
     $$
 where $Y_0$ is the unique solution of equation (4.13). It then
 follows easily that
     $$
      \mbox{sup}_{u\in U} \mbox{inf}_{v \in
V}F(t,x,0,0,0,u, v)=F_0 (t,x,0,0,0)\leq 0,
     $$
and from the definition of $F$\ we see that $W$ is a viscosity
 supersolution of equation (4.1).\\

(2) The second step is devoted to the proof that $W$ is a viscosity
subsolution. For this we suppose that $\varphi \in C^3_{l,b} ([0,T]
\times {\mathbb{R}}^n)$\ and $(t,x)\in [0, T)\times {\mathbb{R}}^n$\
are such that $W-\varphi$\ attains a maximum at $(t,x)$. Without
loss of generality we suppose again $\varphi (t,x)=W(t,x)$. We must
prove that
$$      \mbox{sup}_{u\in U} \mbox{inf}_{v \in
V}F(t,x,0,0,0,u, v)=F_0 (t,x,0,0,0)\geq 0.
     $$
Let us suppose that this is not true. Then there exists some
$\theta>0$ such that \be F_0 (t,x,0,0,0)=\mbox{sup}_{u\in U}
\mbox{inf}_{v \in V}F(t,x,0,0,0,u, v)\leq-\theta<0,\ee and we can
find a measurable function $\psi: U\rightarrow V$ such that
$$ F(t,x,0,0,0,u,\psi(u))\leq -\frac{3}{4}\theta,\ \mbox{for all}\
u\in U.$$ Moreover, since $F(\cdot,x,0,0,0,\cdot,\cdot)$\ is
uniformly continuous on $[0, T]\times U\times V$\ there exists some
$T-t\geq R>0$\ such that \be F(s,x,0,0,0,u,\psi(u))\leq
-\frac{1}{2}\theta,\ \mbox{for all}\ u\in U \mbox{and}\ |s-t|\leq R.
\ee On the other hand, due to the DPP (see Theorem 3.1), for every
$\delta\in(0, R]$,
     $$
     \varphi (t,x) =W(t,x) =\mbox{essinf}_{\beta \in {\mathcal{B}}_{t, t+\delta}}\mbox{esssup}_{u \in
{\mathcal{U}}_{t, t+\delta}}G^{t,x;u,\beta(u)}_{t,t+\delta}
[W(t+\delta, X^{t,x;u,\beta(u)}_{t+\delta})],
     $$
  and from $W\leq\varphi$\ and the monotonicity property of $G^{t,x;u,\beta(u)}_{t,t+\delta}[\cdot]$\ (see Lemma 2.2)\ we obtain
     $$
     \mbox{essinf}_{\beta \in {\mathcal{B}}_{t, t+\delta}}\mbox{esssup}_{u \in
{\mathcal{U}}_{t, t+\delta}} \{G^{t,x;u,\beta(u)}_{t,t+\delta}
[\varphi(t+\delta, X^{t,x;u,\beta(u)}_{t+\delta})] -\varphi
(t,x)\}\geq 0,\ \mbox{P-a.s.}
     $$
 Thus, from Lemma 4.2,
    $$
      \mbox{essinf}_{\beta \in {\mathcal{B}}_{t, t+\delta}}\mbox{esssup}_{u \in
{\mathcal{U}}_{t, t+\delta}} Y^{1,u,\beta(u)}_t \geq 0,\
\mbox{P-a.s.},
    $$
and, in particular,
    $$
     \mbox{esssup}_{u \in
{\mathcal{U}}_{t, t+\delta}} Y^{1,u,\psi(u)}_t \geq 0,\
\mbox{P-a.s.}
    $$
 Here, by putting $\psi_s(u)(\omega)=\psi(u_s(\omega)),\ (s,\omega)\in[t,T]\times
 \Omega$, we identify $\psi$\ as an element of ${\cal{B}}_{t, t+\delta}$.
Given an arbitrarily $\varepsilon>0$ we can choose $u^\varepsilon\in
{\mathcal{U}}_{t, t+\delta}$\ such that
$Y^{1,u^\varepsilon,\psi(u^\varepsilon)}_t\geq -\varepsilon\delta$\
(similar to the proof of (6.21)). From Lemma 4.3 we further have
     \be
   Y^{2,u^\varepsilon,\psi(u^\varepsilon)}_t\geq -C\delta^{\frac{3}{2}}- \varepsilon\delta,\ \mbox{P-a.s.}
     \ee
  Taking into account that $$Y^{2,u^\varepsilon,\psi(u^\varepsilon)}_t
  =E[\int_t^{t+\delta}F(s, x, Y^{2,u^\varepsilon,\psi(u^\varepsilon)}_s,Z^{2,u^\varepsilon,\psi(u^\varepsilon)}_s,
  K^{2,u^\varepsilon,\psi(u^\varepsilon)}_s,u^\varepsilon_s, \psi_s(u^\varepsilon_.))ds|{\mathcal{F}}_{t}]$$\ we get
from the Lipschitz property of $F$ in $(y, z,k)$, (4.20) and Lemma
4.5 that \be\begin{array}{ll}
Y^{2,u^\varepsilon,\psi(u^\varepsilon)}_t&\leq
E[\int_t^{t+\delta}(C|Y^{2,u^\varepsilon,\psi(u^\varepsilon)}_s|+C|Z^{2,u^\varepsilon,\psi(u^\varepsilon)}_s|
+C|\int_EK^{2,u^\varepsilon,\psi(u^\varepsilon)}_s(e)l(x,e)\lambda(de)|\\
&\ \hskip2cm+F(s,
x,0,0,0,u^\varepsilon_s, \psi_s(u^\varepsilon_.)))ds|{\mathcal{F}}_{t}]\\
&\leq C\delta^{\frac{3}{2}}-\frac{1}{2}\theta\delta,\
\mbox{P-a.s.}\end{array}\ee From (4.21) and (4.22),
$-C\delta^{\frac{1}{2}}- \varepsilon\leq
C\delta^{\frac{1}{2}}-\frac{1}{2}\theta,\ \mbox{P-a.s.}$\ Letting
$\delta\downarrow0$,\ and then $\varepsilon\downarrow0$\ we deduce
$\theta\leq0$\ which induces a contradiction.
 Therefore,
     $$
     F_0 (t,x,0,0,0) = \mbox{sup}_{u\in U} \mbox{inf}_{v \in
V}F(t,x,0,0,0,u, v) \geq 0, $$
 and from the definition of $F$, we know that $W$ is a viscosity
 subsolution of equation (4.1). Finally, the results from the first and the second step prove that $W$ is a viscosity
 solution of equation (4.1).\endpf

\br Similarly, we can prove that $U$\ is a viscosity
 solution of equation (4.2). \er

\section{\large Viscosity Solution of Isaacs' Equation: Uniqueness Theorem }
The objective of this section is to study the uniqueness of the
viscosity solution of Isaacs' equation (4.1), \be
 \left \{\begin{array}{ll}
 &\!\!\!\!\! \frac{\partial }{\partial t} W(t,x) +  H^{-}(t, x, W, DW, D^2W)=0,
 \hskip 0.5cm   (t,x)\in [0,T)\times {\mathbb{R}}^n ,  \\
 &\!\!\!\!\!  W(T,x) =\Phi (x), \hskip0.5cm   x \in {\mathbb{R}}^n.
 \end{array}\right.
\ee Recall that
$$\begin{array}{lll}
& H^-(t, x, W, DW, D^2W)= \mbox{sup}_{u \in U}\mbox{inf}_{v \in
V}\{\frac{1}{2}tr(\sigma\sigma^{T}(t, x,
 u, v) D^2W)+ DW.b(t, x, u, v)\\
 &\ \hskip3cm+\int_E(W(t, x+\gamma(t,x,u,v,e))-W(t,x)-DW.\gamma(t,x,u,v,e))\lambda(de)\\
 &+
 f(t, x, W(t,x), DW.\sigma(t,x,u,v),\int_E(W(t, x+\gamma(t,x,u,v,e))-W(t,x))l(x,e)\lambda(de), u,
 v)\},
 \end{array}$$ $\mbox{where}\ t\in [0, T],\ x\in
{\mathbb{R}}^n.$\
 The functions $b, \sigma, f\ \mbox{and}\ \Phi$\ are still supposed to satisfy (H3.1) and (H3.2), respectively.

 We will prove the uniqueness for equation (5.1) in the following
 space of continuous functions

 $\Theta=\{ \varphi\in C([0, T]\times {\mathbb{R}}^n): \exists\ \widetilde{A}>0\ \mbox{such
 that}$ \vskip 0.1cm
 $\mbox{ }\hskip2cm \lim_{|x|\rightarrow \infty}\varphi(t, x)\exp\{-\widetilde{A}[\log((|x|^2+1)^{\frac{1}{2}})]^2\}=0,\
 \mbox{uniformly in}\ t\in [0, T]\}.$  \vskip 0.1cm
\noindent This space of continuous functions is endowed with a
growth condition which is slightly weaker than the assumption of
polynomial growth but more restrictive than that of exponential
growth. This growth condition was introduced by Barles, Buckdahn and
Pardoux~\cite{BBP} and Barles and Imbert~\cite{BI} to prove the uniqueness of the viscosity solution
of an integral-partial differential equation associated with a
decoupled FBSDE with jumps but without controls. It was shown
in~\cite{BBP} that this kind of growth condition is optimal for the
uniqueness and can not be weakened in general. We adapt the ideas
developed in~\cite{BBP} to Isaacs' equation (5.1) to prove the
uniqueness of the viscosity solution in $\Theta$. Since the proof of
the uniqueness in $\Theta$\ for equation (4.2) is essentially the
same we will restrict ourselves to that of (5.1). Before stating the
main result of this section, let us begin with two auxiliary
lemmata. Denoting by $K$\ a Lipschitz constant of $f(t,.,.,.,.,u,
v)$, which is uniformly in $(t,u, v),$\ we have the following

\bl\mbox{  } Let $u_1 \in \Theta$\ be a viscosity subsolution and
$u_2 \in \Theta$\ be a viscosity supersolution of equation (5.1).
Then the function $\omega:= u_1-u_2$\ is a viscosity subsolution of
the equation
    \be
    \left \{
     \begin{array}{lll}
     &\!\!\!\!\!\frac{\partial }{\partial t} \omega(t,x) + \mbox{sup}_{u \in
U, v \in V}\{ A^{u,v}\omega(t,x)+
 B^{u,v}\omega(t,x)+K|\omega(t,x)|+K|D\omega(t,x).\sigma(t, x, u, v)|\\
 &\!\!\!\!\!\mbox{ }\hskip1cm +K(C^{u,v}\omega(t,x))^+ \}= 0, \ \hskip2cm  (t, x)\in [0, T)\times
 {\mathbb{R}}^n,\\
&\!\!\!\!\!\omega(T,x) =0,\ \hskip1cm  x \in {\mathbb{R}}^n.
     \end{array}\right.
   \ee
  \el
\noindent The proof of this lemma follows directly from that of
Lemma 3.7 in~\cite{BBP} with the help of Lemma 1 (Nonlocal Jensen-Ishii's Lemma) in Barles and
Imbert~\cite{BI}.

\noindent Now we can prove the uniqueness theorem.
\bt\mbox{ } We assume that (H3.1) and (H3.2) hold. Let $u_1$ (resp.,
$u_2$) $\in \Theta$
  be a viscosity subsolution (resp., supersolution) of equation
  (5.1). Then we have
     \be
     u_1 (t,x) \leq u_2 (t,x) , \hskip 0.5cm \mbox{for all}\ \ (t,x) \in [0,T] \times {\mathbb{R}}^n .
    \ee
\et \noindent {\bf Proof.} Let us first suppose that $u_1$\ and $u_2$\ are bounded and put $\omega_1:= u_1-u_2$. Theorem 4.1 in~\cite{BI} establishes a
comparison principle for bounded sub- and supersolutions of Hamilton-Jacobi-Bellman equations with
nonlocal term of type (5.2). We know from Lemma 5.1 that $\omega_1$ is a
viscosity subsolution of equation (5.2). On the other hand,
$\omega_2=0$ is, obviously, a viscosity solution and, hence, also a
viscosity supersolution of equation (5.2). Both functions $\omega_1$
and $\omega_2$ are bounded, and the comparison principle stated
in Theorem 4.1 in~\cite{BI} yields that
$u_1-u_2=\omega_1\le\omega_2=0$, i.e., $u_1\le u_2$ on $[0,T]\times
R^n$. Finally, if $u_1,u_2$ are viscosity solutions of (5.2), they
are both viscosity sub- and supersolution, and from the just proved
comparison result we get the equality of $u_1$ and $u_2$. However, under our stand assumptions we can not expect that $W$\ is bounded, so that we have to prove the theorem for $u_1, u_2\in \Theta$. For the proof the following auxiliary lemma is needed.

In analogy to~\cite{BBP} we also have

 \bl\mbox{  }For any
$\widetilde{A}>0,$\ there exists $C_1>0$\ such that the function
$$\chi(t,x)=\exp[(C_1(T-t)+\widetilde{A})\psi(x)],$$
with
$$\psi(x)=[\log((|x|^2+1)^{\frac{1}{2}})+1]^2,\ x\in {\mathbb{R}}^n,$$
satisfies
    \be
     \begin{array}{lll}
     &\frac{\partial }{\partial t}\chi(t,x) + \mbox{sup}_{u \in
U, v \in V}\{ A^{u,v}\chi(t, x)+ B^{u,v}\chi(t, x)+K\chi(t,x) +\\
 & K|D\chi(t,x).\sigma(t, x, u, v)|+K(C^{u,v}\chi(t,x))^+  \}< 0 \ \  \mbox{in}\ [t_1, T]\times
 {\mathbb{R}}^n,\ \mbox{where}\ \ t_1=T-\frac{\widetilde{A}}{C_1}.
     \end{array}
    \ee
   \el
\noindent {\bf Proof.} By direct calculus we first deduce the
following estimates for the first and second derivatives of
$\psi$:
$$ |D\psi(x)|\leq \frac{2[\psi(x)]^{\frac{1}{2}}}{(|x|^2+1)^{\frac{1}{2}}}\leq
4,\ \ \  |D^2\psi(x)|\leq
\frac{C(1+[\psi(x)]^{\frac{1}{2}})}{|x|^2+1},\ \ \ x\in
{\mathbb{R}}^n.$$ These estimates imply that, if $t\in [t_1, T],$
$$
\begin{array}{lll}
     |D\chi(t,x)|&\leq (C_1(T-t)+\widetilde{A})\chi(t,x)|D\psi(x)|\\
 & \leq
 C\chi(t,x)\frac{[\psi(x)]^{\frac{1}{2}}}{(|x|^2+1)^{\frac{1}{2}}},
     \end{array}
$$
and, similarly
$$ |D^2\chi(t,x)| \leq C\chi(t,x)\frac{\psi(x)}{|x|^2+1}.
$$
We should notice that the above estimates do not depend on $C_1$\
because of the definition of $t_1$. Then, since $\gamma$\ is bounded
and since $\psi$\ is Lipschitz continuous in ${\mathbb{R}}^n$, we
have after a long but straight-forward calculus,
$$\chi(t, x+\gamma(t,x,u,v,e))-\chi(t,x)-D\chi(t, x).\gamma(t,x,u,v,e)
\leq C\chi(t,x)\frac{\psi(x)}{|x|^2+1}|\gamma(t,x,u,v,e)|^2,$$ and
$$\chi(t, x+\gamma(t,x,u,v,e))-\chi(t,x)
\leq
C\chi(t,x)\frac{[\psi(x)]^{\frac{1}{2}}}{(|x|^2+1)^{\frac{1}{2}}}|\gamma(t,x,u,v,e)|.$$
In virtue with the above estimates we have
$$  \begin{array}{lll}
     &\frac{\partial }{\partial t}\chi(t,x) + \mbox{sup}_{u \in
U, v \in V}\{ A^{u,v}\chi(t, x)+ B^{u,v}\chi(t, x)+K\chi(t,x) +\\
 & K|D\chi(t,x).\sigma(t, x, u, v)|+K(C^{u,v}\chi(t,x))^+ \}\\
 &\leq -\chi(t,x)\{C_1\psi(x)-C\psi(x)-C[\psi(x)]^{\frac{1}{2}}-\frac{C\psi(x)}{|x|^2+1}-K
 -CK[\psi(x)]^{\frac{1}{2}}-CK\frac{[\psi(x)]^{\frac{1}{2}}}{(|x|^2+1)^{\frac{1}{2}}}\}\\
 &< -\chi(t,x)\{C_1-[C+K]\}\psi(x)< 0,\ \mbox{if}\ C_1>C+K \ \mbox{large
 enough}.
     \end{array}
   $$
\endpf

Now we can continue to prove the uniqueness theorem---Theorem 5.1.\\

\noindent {\bf Proof of Theorem 5.1. (continued)} Let us put $\omega:= u_1-u_2$. Then we
have, for some $\widetilde{A}>0$,
$$\lim_{|x|\rightarrow \infty}\omega(t, x)e^{-\widetilde{A}[\log((|x|^2+1)^{\frac{1}{2}})]^2}=0,$$
uniformly with respect to $t\in [0, T]$. This implies, in
particular, that for any $\alpha >0$, $\omega(t, x)-\alpha\chi(t,
x)$\ is bounded from above in $[t_1, T]\times {\mathbb{R}}^n,$\
and that
$$ M:=\max_{[t_1, T]\times {\mathbb{R}}^n}(\omega-\alpha\chi)(t, x)e^{-K(T-t)}$$
is achieved at some point $(t_0, x_0)\in [t_1, T]\times
{\mathbb{R}}^n$\ (depending on $\alpha$). We now have to distinguish between two cases.\\
For the first case we suppose that: $\omega(t_0, x_0)\leq 0$, for any $\alpha>0$.\\
Then, obviously $M\leq 0$ and $u_1(t, x)-u_2(t, x)\leq
\alpha\chi(t, x)$\ in $[t_1, T]\times {\mathbb{R}}^n$.
Consequently, letting $\alpha$\ tend to zero we obtain
$$u_1(t, x)\leq u_2(t, x),\ \ \mbox{for all}\ (t, x)\in [t_1, T]\times {\mathbb{R}}^n. $$
For the second case we assume that there exists some $\alpha>0$\ such that $\omega(t_0, x_0)> 0$.\\
We notice that $\omega(t, x)-\alpha\chi(t,x)\leq (\omega(t_0,
x_0)-\alpha\chi(t_0, x_0))e^{-K(t-t_0)}\ \ \mbox{in}\ \ [t_1,
T]\times {\mathbb{R}}^n.$\ Then, putting
$$\varphi(t, x)=\alpha\chi(t, x)+(\omega-\alpha\chi)(t_0, x_0)e^{-K(t-t_0)}$$
we get $\omega-\varphi\leq 0=(\omega-\varphi)(t_0, x_0)\ \mbox{in}\
\ [t_1, T]\times {\mathbb{R}}^n.$\ Consequently, since $\omega$\ is
a viscosity subsolution of (5.2) from Lemma 5.1  we have
 $$
     \begin{array}{lll}
     &\frac{\partial }{\partial t}\varphi(t_0, x_0) + \mbox{sup}_{u \in
U, v \in V}\{ A^{u,v}\varphi(t_0, x_0))+ B^{u,v}\varphi(t_0, x_0)+\\
  & K|\varphi(t_0, x_0)| +
 K|D\varphi(t_0, x_0).\sigma(t_0, x_0, u, v)|+K(C^{u,v}\varphi(t_0, x_0))^+ \}\geq 0.
     \end{array}
$$
Moreover, due to our assumption that $\omega(t_0, x_0)>0$\ and since
$\omega(t_0, x_0)=\varphi(t_0, x_0)$\ we can replace $K|\varphi(t_0,
x_0)|$\ by $K\varphi(t_0, x_0)$\ in the above formula. Then, from
the definition of $\varphi$\ and Lemma 5.2,
$$
     \begin{array}{lll}
     &0\leq \alpha\{\frac{\partial \chi}{\partial t} (t_0, x_0) + \mbox{sup}_{u \in
U, v \in V}\{ A^{u,v}\chi(t_0, x_0))+ B^{u,v}\chi(t_0, x_0)+ \\
 & +K\chi(t_0, x_0) + K|D\chi(t_0, x_0).\sigma(t_0, x_0, u, v)|+K(C^{u,v}\chi(t_0, x_0))^+ \}\}< 0
     \end{array}
   $$
which is a contradiction. Finally, by applying successively the
same argument on the interval $[t_2, t_1]$\ with
$t_2=(t_1-\frac{\widetilde{A}}{C_1})^{+},$\ and then, if $t_2>0,$\
on $[t_3, t_2]$\ with $t_3=(t_2-\frac{\widetilde{A}}{C_1})^{+},$\
etc. We get
$$  u_1 (t,x) \leq u_2 (t,x) , \hskip 0.5cm (t,x) \in [0,T] \times {\mathbb{R}}^n .$$
Thus, the proof is complete.\endpf
 \br\mbox{  } Obviously, since the lower value function $W(t,x)$\ is of at most linear growth it belongs to $\Theta$,
   and so $W(t,x)$ is the unique viscosity solution in $\Theta$ of equation (5.1). Similarly we get that the upper value function
   $U(t,x)$\ is the unique viscosity solution in $\Theta$ of equation (4.2).\er

\br\mbox{  } If the Isaacs' condition holds, that is, if for all
$(t, x)\in [0, T]\times {\mathbb{R}}^n,$
$$H^-(t, x, \Psi(t,x), D\Psi(t,x), D^2\Psi(t,x))=H^+(t, x, \Psi(t,x), D\Psi(t,x), D^2\Psi(t,x)),$$
then the equations (5.1) and (4.2) coincide and from the uniqueness
in $\Theta$ of viscosity solution it follows that the lower value
function $W(t,x)$ equals to the upper value function $U(t,x)$ which
means the associated stochastic differential game has a value.\er

\section{\large{Appendix}}
\subsection{FBSDEs with Jumps}
 \hskip1cm In this subsection we give an overview over basic results on BSDEs with jumps associated
 with Forward SDEs with jumps (for short: FBSDEs) for reader's convenience. We consider measurable functions $b: [0,T]\times \Omega\times
{\mathbb{R}}^n \rightarrow {\mathbb{R}}^n \ $,
         $\sigma: [0,T]\times \Omega\times {\mathbb{R}}^n\rightarrow {\mathbb{R}}^{n\times
         d}$\ and $\gamma: [0,T]\times \Omega\times {\mathbb{R}}^n \times E  \rightarrow
         {\mathbb{R}}^n$\
which are supposed to satisfy the following conditions:
 $$
 \begin{array}{ll}
   &\begin{array}{ll}
\mbox{(i)}&b(\cdot,0)\ \mbox{and}\ \sigma(\cdot,0)\ \mbox{are} \
{\cal{F}}_t\mbox{-adapted processes, and there exists some}\\
 & \mbox{constant}\ C>0\  \mbox{such that}\\
 &\hskip 1cm|b(t,0)|+|\sigma(t,0)|\leq C, \ \  \mbox{dtdP-}a.e.;\\
  \end{array}\\
  &\begin{array}{ll}
\mbox{(ii)}&b\ \mbox{and}\ \sigma\ \mbox{are Lipschitz in}\ x,\ \mbox{i.e., there is some constant}\ C>0\ \mbox{such that}\\
           &\hskip 1cm|b(t,x)-b(t,x')|+|\sigma(t,x)-\sigma(t,x')|\leq C| x-x'|,\ \mbox{dtdP-}a.e.,\\
 & \hbox{ \ \ }\hskip9cm\mbox{for }\ x,\ x'\in {\mathbb{R}}^n;\\
\mbox{(iii)}&\mbox{There exists a measurable function}\ \rho :E
\rightarrow {\mathbb{R}}^{+}\mbox{ with }\int_E \rho^2
(e)\lambda(de)<+\infty,\\[.1cm]
   &\mbox{ such that,}\ \ \mbox{for any }\ x, y \in {\mathbb{R}}^n\mbox{ and }\ e \in E,\\[.1cm]
    \end{array}\\
  & \begin{array}{rcl}
 |\gamma (t,x,e)- \gamma (t,y,e)| & \leq & \rho (e)|x-y|, \\[.1cm]
 \ \hskip0.85cm\gamma(\cdot,e)\ \mbox{is} \
{\cal{F}}_t\mbox{-predictable, and}\ |\gamma (t,x,e)| & \leq & \rho
(e)(1+|x|), \ \mbox{dtdP-}a.e..
   \end{array}\\
 \end{array}
  \eqno{\mbox{(H6.1)}}
  $$\par
  We now consider the following SDE with jumps parameterized by the
  initial condition $(t,\zeta)\in[0,T]\times L^2(\Omega,{\cal{F}}_t,P;{\mathbb{R}}^n)$:
  \be
  \left\{
  \begin{array}{rcl}
  dX_s^{t,\zeta}&=&b(s,X_s^{t,\zeta})ds+\sigma(s,X_s^{t,\zeta})dB_s+
                   \int_E \gamma (s,X^{t,\zeta}_{s-},e) \tilde{\mu} (ds,de), \\
  X_t^{t,\zeta}&=&\zeta, \ \hskip6cm\ s\in[t,T].
  \end{array}
  \right.
  \ee
Under the assumption (H6.1), SDE (6.1) has a unique strong solution
and, there exists $C\in {\mathbb{R}}^+$\ such that, for any
$t\in[0,T]\ \mbox{and}\ \zeta,\zeta'\in
L^2(\Omega,{\cal{F}}_t,P;{\mathbb{R}}^n),$
 \be
 \begin{array}{rcl}
 E[\sup\limits_{t\leq s\leq T}| X_s^{t,\zeta}-X_s^{t,\zeta'}|^2|{\cal{F}}_t]
                             &\leq& C|\zeta-\zeta'|^2, \ \ a.s.,\\
  E[\sup\limits_{t\leq s\leq T}| X_s^{t,\zeta}|^2|{\cal{F}}_t]
                       &\leq& C(1+|\zeta|^2),\ \  a.s.
 \end{array}
\ee (Referred to Proposition 1.1 in~\cite{BBP}). We emphasize that
the constant $C$ in (6.2) only depends on the Lipschitz and the
growth constants of $b$, $\sigma$\ and $\gamma$.

 Let now be given two real valued functions $f(t,x,y,z,k)$\ and $\Phi(x)$\ which shall satisfy the
following conditions:
$$
\begin{array}{ll}
\mbox{(i)}&\Phi:\Omega\times {\mathbb{R}}^n\rightarrow {\mathbb{R}}
\ \mbox{is an}\ {\cal{F}}_T\otimes{\cal{B}}({\mathbb{R}}^n)
             \mbox{-measurable random variable and}\\
          & f:[0,T]\times \Omega\times {\mathbb{R}}^n\times {\mathbb{R}}\times
          {\mathbb{R}}^d\times L^2(E, {\cal{B}}(E), \lambda;
{\mathbb{R}}) \rightarrow {\mathbb{R}}\ \mbox{is}\ {\cal{P}}\otimes
{\cal{B}}({\mathbb{R}}^n) \otimes {\cal{B}}({\mathbb{R}})\otimes
{\cal{B}}({\mathbb{R}}^d)\\
& \otimes\ {\cal{B}}(L^2(E, {\cal{B}}(E),
\lambda; {\mathbb{R}}))\mbox{-measurable}.\\
\mbox{(ii)}&\mbox{There exists a constant}\ C>0\ \mbox{such that}\\
          &| f(t,x,y,z,k)-f(t,x',y',z',k')| +| \Phi(x)-\Phi(x')|\\
          &\hskip 4cm\leq C(|x-x'|+ |y-y'|+|z-z'|+||k-k'||),\ \  a.s.,\\
&\hskip 0.5cm\mbox{for all}\ 0\leq t\leq T,\ x,\ x'\in
{\mathbb{R}}^n,\ y,\ y'\in {\mathbb{R}},\ z, \ z'\in {\mathbb{R}}^d\
\mbox{and}\ k,\ k'\in L^2(E, {\cal{B}}(E), \lambda;
{\mathbb{R}}).\\
\mbox{(iii)}&f\ \mbox{and}\ \Phi \ \mbox{satisfy a linear growth condition, i.e., there exists some}\ C>0\\\
    & \mbox{such that, dt}\times \mbox{dP-a.e.},\ \mbox{for all}\ x\in
    {\mathbb{R}}^n,\\
    &\hskip 2cm|f(t,x,0,0,0)| + |\Phi(x)| \leq C(1+|x|).\\
\end{array}
\eqno{\mbox{(H6.2)}} $$
 With the help of the above assumptions we can verify that the coefficient $f(s,X_s^{t,\zeta},y,z,k)$\
 satisfies the hypotheses (A1) and (A2) and  $\xi=\Phi(X_T^{t,\zeta})$ $\in
 L^2(\Omega,{\cal{F}}_T,P;{\mathbb{R}})$. Therefore, the following
 BSDE with jump possesses a unique solution:
\be \left\{
\begin{array}{rcl}
-dY_s^{t,\zeta}&=&f(s,X_s^{t,\zeta},Y_s^{t,\zeta},Z_s^{t,\zeta},
K_s^{t,\zeta})ds-Z_s^{t,\zeta}dB_s-\int_E K_{s-}^{t,\zeta}(e) \tilde{\mu}(ds,de),\\
Y_T^{t,\zeta}&=&\Phi(X_T^{t,\zeta}),\ \hskip5cm s\in [t, T].\\
\end{array}
\right. \ee

\bp  We suppose that the hypotheses (H6.1) and (H6.2) hold. Then,
for any $0\leq t\leq T$ and the associated initial conditions
 $\zeta,\zeta'\in L^2(\Omega,{\cal{F}}_t,P;{\mathbb{R}}^n)$, we
have the following estimates:\\  $\mbox{}\hskip3cm\mbox{\rm(i)}
E[\sup\limits_{t\leq s\leq T}|Y_s^{t,\zeta}|^2 +
\int_t^T|Z_s^{t,\zeta}|^2ds|{{\cal{F}}_t}]+E[
\int_t^T\int_E|K_s^{t,\zeta}(e)|^2\lambda(de)ds|{{\cal{F}}_t}]$\\
$\mbox{  } \hskip4cm\leq C(1+|\zeta|^2),\
a.s.; $\\
$\mbox{}\hskip3cm\mbox{\rm(ii)}E[\sup\limits_{t\leq s\leq
T}|Y_s^{t,\zeta}-Y_s^{t,\zeta'}|^2+\int_t^T|Z_s^{t,\zeta}-Z_s^{t,\zeta'}|^2ds|{{\cal{F}}_t}]+\\
\mbox{}\hskip4cm
 E[\int_t^T\int_E|K_s^{t,\zeta}(e)-K_s^{t,\zeta'}(e)|^2\lambda(de)ds|{{\cal{F}}_t}]\leq
C|\zeta-\zeta'|^2,\  a.s. $\\ In particular, \be
 \begin{array}{lll}
\mbox{\rm(iii)}&|Y_t^{t,\zeta}|\leq C(1+|\zeta|),\  a.s.; \hskip3cm\\
\mbox{\rm(iv)}&|Y_t^{t,\zeta}-Y_t^{t,\zeta'}|\leq C|\zeta-\zeta'|,\  a.s.\hskip3cm \\
\end{array}
\ee The above constant $C>0$\ depends only on the Lipschitz and the
growth constants of $b$,\ $\sigma$, $\gamma$, $f$\ and $\Phi$. \ep
 \noindent From Lemma 2.1, the estimate (6.2) and It\^{o}'s formula we can prove this proposition . \vskip 0.3cm
 Let us now introduce the random field:
\be u(t,x)=Y_s^{t,x}|_{s=t},\ (t, x)\in [0, T]\times{\mathbb{R}}^n,
\ee where $Y^{t,x}$ is the solution of BSDE (6.3) with $x \in
{\mathbb{R}}^n$\ at the place of $\zeta\in
L^2(\Omega,{\cal{F}}_t,P;{\mathbb{R}}^n).$\\

As a consequence of Proposition 6.1 we have that, for all $t \in [0,
T] $, P-a.s.,
 \be
\begin{array}{ll}
\mbox{(i)}&| u(t,x)-u(t,y)| \leq C|x-y|,\ \mbox{for all}\ x, y\in {\mathbb{R}}^n;\\
\mbox{(ii)}&| u(t,x)|\leq C(1+|x|),\ \mbox{for all}\ x\in {\mathbb{R}}^n.\\
\end{array}
\ee
 \br\ In the general situation $u$\ is random and forms an adapted random field, that is, for any $x\in
 {\mathbb{R}}^n,\ u(\cdot,x)$ is an ${\cal{F}}_t-$adapted real valued process. Indeed, recall that
  $b, \sigma, f\ \mbox{and}\ \Phi$ all are ${\cal{F}}_t$-adapted random functions. On the other hand, it is well
  known that, under the additional assumption that the
  functions $$ b, \sigma, \gamma, f\ \mbox{and}\ \Phi\ \mbox{are
  deterministic,}\eqno{\mbox{(H6.3)}}  $$
 $u$ is also a deterministic function of $(t,x)$. \er

 The random field $u$\ and $Y^{t,\zeta},\ (t, \zeta)\in [0,
T]\times L^2(\Omega,{\cal{F}}_t,P;{\mathbb{R}}^n),$\ are related by
the following theorem.
 \bt Under the assumptions (H6.1) and (H6.2), for any $t\in [0, T]$\ and $\zeta\in
L^2(\Omega,{\cal{F}}_t,P;{\mathbb{R}}^n),$\ we have \be
u(t,\zeta)=Y_t^{t,\zeta},\ \mbox{ P-a.s.}. \ee \et

\br \hskip0.1cm Obviously,\ $ Y_s^{t, \zeta}=Y_s^{s,
X_s^{t,\zeta}}=u(s, X_s^{t,\zeta}),\ \mbox{ P-a.s.}.$ \er

The proof of Theorem 6.1 is similar to the proof in Peng~\cite{P1}
for the FBSDE with Brownian motion, also can refer to Theorem A.1 in~\cite{LP}, for reader's convenience we give the proof here. It makes use of the following
definition.

 \noindent
\bde For any t $\in [0, T]$, a sequence $\{A_i\}_{i=1}^{N}\subset
{\cal{F}}_t\ (\mbox{with}\ 1\leq N\leq \infty)$ is called a
partition of $(\Omega, {\cal{F}}_t)$\ if\ \
$\cup_{i=1}^{N}A_i=\Omega$\ and $ A_i\cap A_j=\phi, \
\mbox{whenever}\ i\neq j.$ \ede
 \noindent \textbf{Proof} (of Theorem 6.1): We first consider
the case where $\zeta$ is a simple random variable of the form\
\be\zeta=\sum\limits^N\limits_{i=1}x_i\textbf{1}_{A_i},
                        \ee
where$\{A_i\}^N_{i=1}$\ is a finite partition of $(\Omega,{\cal{F}}_t)$\ and $x_i\in {\mathbb{R}}^n$,\ for $1\leq i\leq N.$\\
For each $i$, we put $(X_s^i,Y_s^i,Z_s^i,K_s^i)\equiv
                 (X_s^{t,x_i},Y_s^{t,x_i},Z_s^{t,x_i},K_s^{t,x_i}).$ Then $X^i$ is the solution of the SDE
$$
X^i_s =x_i +\int^s_t b(r,X^i_r)dr +\int^s_t \sigma
(r,X^i_r)dB_r+\int^s_t\int_E\gamma(r,X^i_{r-},e)\widetilde{\mu}(dr,de),\
s\in [t,T],
$$\\
 and $(Y^i,Z^i,K^i)$ is the solution of the associated BSDE
$$
Y^i_s =\Phi(X^i_T) +\int^T_s f(r,X^i_r,Y^i_r,Z^i_r,K^i_r)dr
   -\int^T_s Z^i_r dB_r-\int^T_s\int_EK^i_r(e)\widetilde{\mu}(dr,de),\ s\in [t,T].
$$
The above two equations are multiplied by $\textbf{1}_{A_i}$\ and
summed up with respect to $i$. Thus, taking into account that
$\sum\limits_i \varphi (x_i)\textbf{1}_{A_i}=\varphi (\sum\limits_i
x_i \textbf{1}_{A_i})$, we get
$$
\begin{array}{rcl}
 \sum\limits_{i=1}\limits^{N} \textbf{1}_{A_i} X^i_s &=&\zeta+ \int^s_t b(r,\sum\limits _{i=1}\limits^{N} \textbf{1}_{A_i} X^i_r )dr
  +\int^s_t \sigma (r,\sum\limits_{i=1}\limits^{N} \textbf{1}_{A_i}
  X^i_r)dB_r\\
  &&+\int^s_t\int_E\gamma(r,\sum\limits _{i=1}\limits^{N} \textbf{1}_{A_i}
  X^i_{r-},e)\widetilde{\mu}(dr,de);
\end{array}
$$and
$$
\begin{array}{rcl}
\sum\limits _{i=1}\limits^{N}\textbf{1}_{A_i} Y^i_s & = &
\Phi(\sum\limits _{i=1}\limits^{N} \textbf{1}_{A_i} X^i_T)+\int^T_s
f(r,\sum \limits_{i=1}^{N} \textbf{1}_{A_i} X^i_r,
  \sum\limits _{i=1}^{N} \textbf{1}_{A_i} Y^i_r,
 \sum\limits _{i=1}^{N} \textbf{1}_{A_i} Z^i_r,\sum\limits _{i=1}^{N} \textbf{1}_{A_i} K^i_r)dr \\
  & & -\int^T_s \sum\limits _{i=1}^{N} \textbf{1}_{A_i} Z^i_r
  dB_r-\int^T_s\int_E \sum\limits _{i=1}^{N} \textbf{1}_{A_i} K^i_r(e)\widetilde{\mu}(dr,de).
\end{array}
$$
Then the strong uniqueness property of the solution of the SDE and
the BSDE yields
$$
X^{t,\zeta}_s =\sum \limits_{i=1}^{N} X^i_s \textbf{1}_{A_i},\
(Y^{t,\zeta}_s,Z^{t,\zeta}_s,K^{t,\zeta}_s) =(\sum \limits_{i=1}^{N}
\textbf{1}_{A_i} Y^i_s, \sum \limits_{i=1}^{N} \textbf{1}_{A_i}
Z^i_s,\sum \limits_{i=1}^{N} \textbf{1}_{A_i} K^i_s),\ s\in [t, T].
$$
Finally, from $u(t,x_i)=Y^i_t,\ 1\leq i\leq N$, we deduce that
$$
Y^{t,\zeta}_t=\sum \limits_{i=1}^{N}
Y^i_t\textbf{1}_{A_i}=\sum\limits_{i=1}^{N}u(t,x_i) \textbf{1}_{A_i}
=u(t,\sum \limits_{i=1}^{N} x_i \textbf{1}_{A_i}) =u(t,\zeta).
$$
Therefore, for simple random variables, we have the desired result.

Given a general $\zeta\in L^2 (\Omega ,{\mathcal{F}}_t
,P;{\mathbb{R}}^n)$ we can choose a sequence of simple random
variables $\{\zeta_i\}$ which
 converges to $\zeta$ in $L^2(\Omega ,{\mathcal{F}}_t
,P;{\mathbb{R}}^n)$. Consequently, from the estimates (6.4), (6.6)
and the first step of the proof, we have
$$
\begin{array}{lrcl}
&E|Y^{t,\zeta_i}_t-Y^{t,\zeta}_t|^2&\leq&CE|\zeta_i -\zeta|^2\rightarrow 0,\ i\rightarrow\infty,\\
\mbox{ }\hskip1cm&
E|u(t,\zeta_i)-u(t,\zeta)|^2 &\leq& CE|\zeta_i -\zeta|^2 \rightarrow 0,\ i\rightarrow\infty,\\
\hbox{and}\hskip1cm& Y^{t,\zeta_i}_t&=& u(t,\zeta_i),\ i\geq 1.
\end{array}
$$
Then the proof is complete.\endpf

\subsection{The Proof of Theorem 3.1}

\noindent \textbf{Proof.}  To simplify notations we put
$$W_\delta(t,x) =\mbox{essinf}_{\beta \in {\mathcal{B}}_{t,
t+\delta}}\mbox{esssup}_{u \in {\mathcal{U}}_{t,
t+\delta}}G^{t,x;u,\beta(u)}_{t,t+\delta} [W(t+\delta,
X^{t,x;u,\beta(u)}_{t+\delta})].$$ In analogy to $W(t,x)$\ it can be
easily shown that $W_\delta(t,x)$\ is well-defined. The proof that
$W_\delta(t,x)$\ coincides with $W(t,x)$\ will be split into a
sequel of lemmas which all supposed that (H3.1) and (H3.2) are
satisfied.

\bl $W_\delta(t,x)$\ is deterministic.\el The proof of this lemma
uses the same ideas as that of Proposition 3.1 so  it is
omitted here.\endpf

\bl$W_\delta(t,x)\leq W(t,x).$\el

\noindent\textbf{Proof}. Let $\beta\in {\mathcal{B}}_{t, T}$\ be
arbitrarily fixed. Then, given a $u_2(\cdot)\in
{\mathcal{U}}_{t+\delta, T},$\ we define as follows the restriction
$\beta_1$\ of $\beta$\ to ${\mathcal{U}}_{t+\delta, T}:$
$$\beta_1(u_1):=\beta(u_1\oplus u_2 )|_{[t,
t+\delta]},\ \mbox{ }\ u_1(\cdot)\in {\mathcal{U}}_{t, t+\delta},
$$
where $u_1\oplus u_2:=u_1\textbf{1}_{[t,
t+\delta]}+u_2\textbf{1}_{(t+\delta, T]}$\ extends $u_1(\cdot)$\
to an element of ${\mathcal{U}}_{t, T}$. It is easy to check that
$\beta_1\in {\mathcal{B}}_{t, t+\delta}.$\ Moreover, from the
nonanticipativity property of $\beta$\ we deduce that $\beta_1$\
is independent of the special choice of $u_2(\cdot)\in
{\mathcal{U}}_{t+\delta, T}.$\ Consequently, from the definition
of $W_\delta(t,x),$
 \be W_\delta(t,x)\leq \mbox{esssup}_{u_1
\in {\mathcal{U}}_{t,
t+\delta}}G^{t,x;u_1,\beta_1(u_1)}_{t,t+\delta} [W(t+\delta,
X^{t,x;u_1,\beta_1(u_1)}_{t+\delta})],\ \mbox{P-a.s.} \ee We use the
notation $I_\delta(t, x, u, v):=G^{t,x;u,v}_{t,t+\delta}
[W(t+\delta, X^{t,x;u,v}_{t+\delta})]$\ and notice that there exists
a sequence $\{u_i^1,\ i\geq 1\}\subset {\mathcal{U}}_{t, t+\delta}$\
such that
$$I_\delta(t, x, \beta_1):=\mbox{esssup}_{u_1 \in {\mathcal{U}}_{t,
t+\delta}}I_\delta(t, x, u_1, \beta_1(u_1))=\mbox{sup}_{i\geq
1}I_\delta(t, x, u_i^1, \beta_1(u_i^1)),\ \ \mbox{P-a.s.}.$$ For any
$\varepsilon>0,$\ we put $\widetilde{\Gamma}_i:=\{I_\delta(t, x,
\beta_1)\leq I_\delta(t, x, u_i^1, \beta_1(u_i^1))+\varepsilon\}\in
{\mathcal{F}}_{t},\ i\geq 1.$\ Then
$\Gamma_1:=\widetilde{\Gamma}_1,\
\Gamma_i:=\widetilde{\Gamma}_i\backslash(\cup^{i-1}_{l=1}\widetilde{\Gamma}_l)\in
{\mathcal{F}}_{t},\ i\geq 2,$\ form an $(\Omega,
{\mathcal{F}}_{t})$-partition, and $u^\varepsilon_1:=\sum_{i\geq
1}\textbf{1}_{\Gamma_i}u_i^1$\ belongs obviously to
${\mathcal{U}}_{t, t+\delta}.$\ Moreover, from the nonanticipativity
of $\beta_1$\ we have $\beta_1(u^\varepsilon_1)=\sum_{i\geq
1}\textbf{1}_{\Gamma_i}\beta_1(u_i^1),$\ and from the uniqueness of
the solution of the FBSDE, we deduce that $I_\delta(t, x,
u^\varepsilon_1, \beta_1(u^\varepsilon_1))=\sum_{i\geq
1}\textbf{1}_{\Gamma_i}I_\delta(t, x, u_i^1, \beta_1(u_i^1)),\
\mbox{P-a.s.}$\ Hence, \be
\begin{array}{llll}
W_\delta(t,x)\leq I_\delta(t, x, \beta_1)&\leq &\sum_{i\geq
1}\textbf{1}_{\Gamma_i}I_\delta(t, x, u_i^1, \beta_1(u_i^1))
+\varepsilon=I_\delta(t, x, u^\varepsilon_1,
\beta_1(u^\varepsilon_1))+\varepsilon\\
&=& G^{t,x;u^\varepsilon_1, \beta_1(u^\varepsilon_1)}_{t,t+\delta}
[W(t+\delta, X^{t,x;u^\varepsilon_1,
\beta_1(u^\varepsilon_1)}_{t+\delta})]+\varepsilon,\ \mbox{P-a.s.}
\end{array}
\ee
 On the other hand, using the fact that $\beta_1(\cdot):=\beta(\cdot\oplus u_2
)\in {\mathcal{B}}_{t, t+\delta}$\ does not depend on $u_2(\cdot)\in
{\mathcal{U}}_{t+\delta, T}$\ we can define
$\beta_2(u_2):=\beta(u^\varepsilon_1\oplus u_2)|_{[t+\delta, T]},\
\mbox{for all }\ u_2(\cdot)\in {\mathcal{U}}_{t+\delta, T}. $\ The
such defined $\beta_2: {\mathcal{U}}_{t+\delta, T}\rightarrow
{\mathcal{V}}_{t+\delta, T}$\ belongs to ${\mathcal{B}}_{t+\delta,
T}\ \mbox{since}\ \beta\in {\mathcal{B}}_{t, T}$. Therefore, from
the definition of $W(t+\delta,y)$\ we have, for any $y\in
{\mathbb{R}}^n,$
$$W(t+\delta,y)\leq \mbox{esssup}_{u_2 \in {\mathcal{U}}_{t+\delta, T}}J(t+\delta, y; u_2, \beta_2(u_2)),\ \mbox{P-a.s.}$$
Finally, because there exists a constant $C\in {\mathbb{R}}$\ such
that \be
\begin{array}{llll}
{\rm(i)} & |W(t+\delta,y)-W(t+\delta,y')| \leq C|y-y'|,\ \mbox{for any}\ y,\ y' \in {\mathbb{R}}^n;  \\
{\rm(ii)} & |J(t+\delta, y, u_2, \beta_2(u_2))-J(t+\delta, y',
u_2, \beta_2(u_2))| \leq C|y-y'|,\ \mbox{P-a.s.,}\\
 &\mbox{ }\hskip1cm \mbox{for any}\ u_2\in {\mathcal{U}}_{t+\delta, T},
\end{array}
\ee (see Lemma 3.3-(i) and (3.6)-(i)) we can show by approximating
$X^{t,x;u_1^\varepsilon,\beta_1(u_1^\varepsilon)}_{t+\delta}$\ that
$$W(t+\delta, X^{t,x;u_1^\varepsilon,\beta_1(u_1^\varepsilon)}_{t+\delta} )\leq
\mbox{esssup}_{u_2 \in {\mathcal{U}}_{t+\delta, T}}J(t+\delta,
X^{t,x;u_1^\varepsilon,\beta_1(u_1^\varepsilon)}_{t+\delta}; u_2,
\beta_2(u_2)),\ \mbox{P-a.s.}$$ To estimate the right side of the
latter inequality we note that there exists some sequence $\{u_j^2,\
j\geq 1\}\subset {\mathcal{U}}_{t+\delta, T}$\ such that
$$\mbox{esssup}_{u_2 \in {\mathcal{U}}_{t+\delta,
T}}J(t+\delta,X^{t,x;u_1^\varepsilon,\beta_1(u_1^\varepsilon)}_{t+\delta};
u_2, \beta_2(u_2))=\mbox{sup}_{j\geq
1}J(t+\delta,X^{t,x;u_1^\varepsilon,\beta_1(u_1^\varepsilon)}_{t+\delta};
u^2_j, \beta_2(u^2_j)),\ \mbox{P-a.s.}$$
 Then, putting\\
$\widetilde{\Delta}_j:=\{\mbox{esssup}_{u_2 \in
{\mathcal{U}}_{t+\delta,
T}}J(t+\delta,X^{t,x;u_1^\varepsilon,\beta_1(u_1^\varepsilon)}_{t+\delta};
u_2, \beta_2(u_2))\leq
J(t+\delta,X^{t,x;u_1^\varepsilon,\beta_1(u_1^\varepsilon)}_{t+\delta};
u^2_j, \beta_2(u^2_j))+\varepsilon\}\in {\mathcal{F}}_{t+\delta},\
j\geq 1;$\ we have with $\Delta_1:=\widetilde{\Delta}_1,\
\Delta_j:=\widetilde{\Delta}_j\backslash(\cup^{j-1}_{l=1}\widetilde{\Delta}_l)\in
{\mathcal{F}}_{t+\delta},\ j\geq 2,$\ an $(\Omega,
{\mathcal{F}}_{t+\delta})$-partition and
$u^\varepsilon_2:=\sum_{j\geq 1}\textbf{1}_{\Delta_j}u_j^2$\
 $\in {\mathcal{U}}_{t+\delta, T}.$ From
the nonanticipativity of $\beta_2$\ we have
$\beta_2(u^\varepsilon_2)=\sum_{j\geq
1}\textbf{1}_{\Delta_j}\beta_2(u_j^2),$\ and from the definition of
$\beta_1,\ \beta_2$\ we know that $\beta(u_1^\varepsilon\oplus
u_2^\varepsilon)=\beta_1(u_1^\varepsilon)\oplus
\beta_2(u_2^\varepsilon ).$\ Thus, again from the uniqueness of the
solution of our FBSDE, we get
$$\begin{array}{lcl}
J(t+\delta,X^{t,x;u_1^\varepsilon,\beta_1(u_1^\varepsilon)}_{t+\delta};
u_2^\varepsilon,
\beta_2(u_2^\varepsilon))&=&Y_{t+\delta}^{t+\delta,X^{t,x;u_1^\varepsilon,\beta_1(u_1^\varepsilon)}_{t+\delta};
u_2^\varepsilon, \beta_2(u_2^\varepsilon)}\ \hskip2cm \mbox{(see (3.8))}\\
&=&\sum_{j\geq
1}\textbf{1}_{\Delta_j}Y_{t+\delta}^{t+\delta,X^{t,x;u_1^\varepsilon,\beta_1(u_1^\varepsilon)}_{t+\delta};
 u_j^2, \beta_2( u_j^2)}\\
&=&\sum_{j\geq
1}\textbf{1}_{\Delta_j}J(t+\delta,X^{t,x;u_1^\varepsilon,\beta_1(u_1^\varepsilon)}_{t+\delta};
u_j^2, \beta_2(u_j^2)),\ \mbox{P-a.s.}
\end{array}$$
Consequently, \be
\begin{array}{lll}
W(t+\delta,
X^{t,x;u_1^\varepsilon,\beta_1(u_1^\varepsilon)}_{t+\delta} )&\leq
& \mbox{esssup}_{u_2 \in {\mathcal{U}}_{t+\delta,
T}}J(t+\delta,X^{t,x;u_1^\varepsilon,\beta_1(u_1^\varepsilon)}_{t+\delta};
u_2, \beta_2(u_2))\\
&\leq& \sum_{j\geq
1}\textbf{1}_{\Delta_j}Y_{t+\delta}^{t,x;u_1^\varepsilon\oplus
u_j^2,\beta(u_1^\varepsilon\oplus
u_j^2)}+\varepsilon\\
& = & Y_{t+\delta}^{t,x;u_1^\varepsilon\oplus u^\varepsilon_2,
\beta(u_1^\varepsilon\oplus
u^\varepsilon_2)}+\varepsilon\\
&=&Y_{t+\delta}^{t,x;u^\varepsilon,\beta(u^\varepsilon)}+\varepsilon,\
 \mbox{P-a.s.,}
\end{array}
\ee where $u^\varepsilon:= u_1^\varepsilon\oplus u^\varepsilon_2\in
{\mathcal{U}}_{t, T}.$\ From (6.10), (6.12), Lemma 2.2 (comparison
theorem for BSDEs) and Lemma 2.3 we have, for some constant $C\in
{\mathbb{R}},$
 \be
\begin{array}{lll}
W_\delta(t,x)&\leq& G^{t,x;u^\varepsilon_1,
\beta_1(u^\varepsilon_1)}_{t,t+\delta}
[Y_{t+\delta}^{t,x;u^\varepsilon,\beta(u^\varepsilon)}+\varepsilon]+\varepsilon \\
&\leq& G^{t,x;u^\varepsilon_1,
\beta_1(u^\varepsilon_1)}_{t,t+\delta}
[Y_{t+\delta}^{t,x;u^\varepsilon,\beta(u^\varepsilon)}]+
(C+1)\varepsilon\\
& =& G^{t,x;u^\varepsilon, \beta(u^\varepsilon)}_{t,t+\delta}
[Y_{t+\delta}^{t,x;u^\varepsilon,\beta(u^\varepsilon)}]+
(C+1)\varepsilon\\
& =& Y_{t}^{t,x;u^\varepsilon,\beta(u^\varepsilon)}+
(C+1)\varepsilon\\
&\leq& \mbox{esssup}_{u \in {\mathcal{U}}_{t,
T}}Y_{t}^{t,x;u,\beta(u)}+ (C+1)\varepsilon,\ \mbox{P-a.s.}
\end{array}
\ee Since $\beta\in {\mathcal{B}}_{t, T}$\ has been arbitrarily
chosen we have (6.13) for all $\beta\in {\mathcal{B}}_{t, T}$.
Therefore,
 \be W_\delta(t,x)\leq \mbox{essinf}_{\beta\in
{\mathcal{B}}_{t, T}}\mbox{esssup}_{u \in {\mathcal{U}}_{t,
T}}Y_{t}^{t,x;u,\beta(u)}+ (C+1)\varepsilon= W(t, x)+
(C+1)\varepsilon.\ee Finally, letting $\varepsilon\downarrow0,\
\mbox{we get}\ W_\delta(t,x)\leq W(t, x).$\endpf

\bl$ W(t, x)\leq W_\delta(t,x).$\el

 \noindent \textbf{Proof}. We continue to use the notations introduced above. From the definition of
$W_\delta(t,x)$\ we have
$$
\begin{array}{lll}
W_\delta(t,x)&=& \mbox{essinf}_{\beta_1 \in {\mathcal{B}}_{t,
t+\delta}}\mbox{esssup}_{u_1 \in {\mathcal{U}}_{t,
t+\delta}}G^{t,x;u_1,\beta_1(u_1)}_{t,t+\delta} [W(t+\delta,
X^{t,x;u_1,\beta_1(u_1)}_{t+\delta})]\\
&=&\mbox{essinf}_{\beta_1 \in {\mathcal{B}}_{t,
t+\delta}}I_\delta(t, x, \beta_1),
\end{array}
$$
and, for some sequence $\{\beta_i^1,\ i\geq 1\}\subset
{\mathcal{B}}_{t, t+\delta},$
$$W_\delta(t,x)=\mbox{inf}_{i\geq
1}I_\delta(t, x, \beta_i^1),\ \mbox{P-a.s.}$$ For any
$\varepsilon>0,$\ we let $\widetilde{\Lambda}_i:=\{I_\delta(t, x,
\beta_i^1)-\varepsilon\leq W_\delta(t,x)\}\in {\mathcal{F}}_{t},\
i\geq 1,$ $\Lambda_1:=\widetilde{\Lambda}_1\ \mbox{and}\
\Lambda_i:=\widetilde{\Lambda}_i\backslash(\cup^{i-1}_{l=1}\widetilde{\Lambda}_l)\in
{\mathcal{F}}_{t},\ i\geq 2.$\ Then $\{\Lambda_i,\ i\geq 1\}$\ is an
$(\Omega, {\mathcal{F}}_{t})$-partition,
$\beta^\varepsilon_1:=\sum_{i\geq
1}\textbf{1}_{\Lambda_i}\beta_i^1$\ belongs to ${\mathcal{B}}_{t,
t+\delta},$\ and from the uniqueness of the solution of our FBSDE we
conclude that $I_\delta(t, x, u_1,
\beta^\varepsilon_1(u_1))=\sum_{i\geq
1}\textbf{1}_{\Lambda_i}I_\delta(t, x, u_1, \beta_i^1(u_1)),$ $\
\mbox{P-a.s., for all}$\ \ $u_1(\cdot)\in {\mathcal{U}}_{t,
t+\delta}.$\ Hence,
 \be
\begin{array}{lll}
W_\delta(t,x)&\geq &\sum_{i\geq
1}\textbf{1}_{\Lambda_i}I_\delta(t, x,
\beta_i^1)-\varepsilon\\
&\geq&\sum_{i\geq 1}\textbf{1}_{\Lambda_i}I_\delta(t, x, u_1,
\beta_i^1(u_1))-\varepsilon \\
&=& I_\delta(t, x, u_1, \beta^\varepsilon_1(u_1))-\varepsilon\\
&=& G^{t,x;u_1, \beta^\varepsilon_1(u_1)}_{t,t+\delta}
[W(t+\delta, X^{t,x;u_1,
\beta_1^\varepsilon(u_1)}_{t+\delta})]-\varepsilon,\ \mbox{P-a.s.,
for all}\ \ u_1\in {\mathcal{U}}_{t, t+\delta}.
\end{array}
\ee
 On the other hand, from the definition of $W(t+\delta,y),$\
with the same technique as before, we deduce that, for any $y\in
{\mathbb{R}}^n,$\ there exists $\beta^\varepsilon_y\in
{\mathcal{B}}_{t+\delta, T}$\ \ such that \be W(t+\delta,y)\geq
\mbox{esssup}_{u_2 \in {\mathcal{U}}_{t+\delta,T}}J(t+\delta, y;
u_2, \beta^\varepsilon_y(u_2))-\varepsilon,\ \mbox{P-a.s.}\ee Let
$\{O_i\}_{i\geq1}\subset {\mathcal{B}}({\mathbb{R}}^n)$\ be a
decomposition of ${\mathbb{R}}^n$\ such that
$\sum\limits_{i\geq1}O_i={\mathbb{R}}^n\ \mbox{and}\
\mbox{diam}(O_i)\leq \varepsilon,\ i\geq 1.$\ Moreover, we fix
arbitrarily for each $i\geq 1$\ an element $y_i$\ of $O_i,\
i\geq1.$\ Then, defining $[X^{t,x;u_1,
\beta_1^\varepsilon(u_1)}_{t+\delta}]:=\sum\limits_{i\geq1}y_i\textbf{1}_{\{X^{t,x;u_1,
\beta_1^\varepsilon(u_1)}_{t+\delta}\in O_i\}},$\ we have \be
|X^{t,x;u_1, \beta_1^\varepsilon(u_1)}_{t+\delta}-[X^{t,x;u_1,
\beta_1^\varepsilon(u_1)}_{t+\delta}]|\leq \varepsilon,\
\mbox{everywhere on}\ \Omega, \ \mbox{for all}\ u_1\in
{\mathcal{U}}_{t, t+\delta}.\ee\ Furthermore, as we have seen above,
for each $y_i$\ there exists some $\beta^\varepsilon_{y_i}\in
{\mathcal{B}}_{t+\delta, T}$\ such that (6.16) holds, and, clearly,
$\beta^{\varepsilon}_{u_1}:=\sum\limits_{i\geq1}\textbf{1}_{\{X^{t,x;u_1,
\beta_1^\varepsilon(u_1)}_{t+\delta}\in
O_i\}}\beta^\varepsilon_{y_i}\in {\mathcal{B}}_{t+\delta, T}.$

 Now we can define the new strategy
$\beta^{\varepsilon}(u):=\beta_1^\varepsilon(u_1)\oplus
\beta^{\varepsilon}_{u_1}(u_2),\ u\in {\mathcal{U}}_{t, T},\
\mbox{where}\ u_1=u|_{[t, t+\delta]},\ u_2=u|_{(t+\delta, T]}$\
(restriction of $u$ to $[t, t+\delta]\times \Omega$\ and $(t+\delta,
T]\times \Omega$, resp.). Obviously, $\beta^{\varepsilon}$\ maps
${\mathcal{U}}_{t,T}$\ into ${\mathcal{V}}_{t,T}.$\ Moreover,
$\beta^{\varepsilon}$\ is nonanticipating: Indeed, let $S:
\Omega\longrightarrow[t, T]$\ be an ${\mathcal{F}}_r$-stopping time
and $u, u'\in {\mathcal{U}}_{t,T}$\ be such that $u\equiv u'$\ on
$\textbf{[\![}t, S\textbf{]\!]}$. Decomposing $u,\ u'$\ into $u_1,
u'_1\in {\mathcal{U}}_{t,t+\delta},\ u_2, u'_2\in
{\mathcal{U}}_{t+\delta, T}$\ such that $u=u_1\oplus u_2\
\mbox{and}\ u'=u'_1\oplus u'_2,$\ we have $u_1\equiv u_1'$\ on
$\textbf{[\![}t, S\wedge(t+\delta)\textbf{]\!]}$\ from which we get
$\beta_1^\varepsilon(u_1)\equiv \beta_1^\varepsilon(u_1')$\ on
$\textbf{[\![}t, S\wedge(t+\delta)\textbf{]\!]}$\ (recall that
$\beta_1^\varepsilon$\ is nonanticipating). On the other hand,
$u_2\equiv u_2'$\ on $\textbf{]\!]}t+\delta,
S\vee(t+\delta)\textbf{]\!]}(\subset (t+\delta,T]\times
\{S>t+\delta\}),$\ and on $\{S>t+\delta\}$\ we have $X^{t,x;u_1,
\beta_1^\varepsilon(u_1)}_{t+\delta}=X^{t,x;u'_1,
\beta_1^\varepsilon(u'_1)}_{t+\delta}.$\ Consequently, from our
definition, $\beta^{\varepsilon}_{u_1}=\beta^{\varepsilon}_{u'_1}$\
on $\{S>t+\delta\}$\ and
$\beta^{\varepsilon}_{u_1}(u_2)\equiv\beta^{\varepsilon}_{u'_1}(u'_2)$\
on $\textbf{]\!]}t+\delta, S\vee(t+\delta)\textbf{]\!]}.$ This
yields $\beta^{\varepsilon}(u)=\beta_1^\varepsilon(u_1)\oplus
\beta^{\varepsilon}_{u_1}(u_2)\equiv\beta_1^\varepsilon(u'_1)\oplus
\beta^{\varepsilon}_{u'_1}(u'_2)=\beta^{\varepsilon}(u')$\ on
$\textbf{[\![}t, S\textbf{]\!]}$, from where it follows that
$\beta^{\varepsilon}\in {\mathcal{B}}_{t, T}.$

Let now $u\in {\mathcal{U}}_{t, T}$\ be arbitrarily chosen and
decomposed into $u_1=u|_{[t, t+\delta]}\in {\mathcal{U}}_{t,
t+\delta}$\ and $u_2=u|_{(t+\delta, T]}\in {\mathcal{U}}_{t+\delta,
T}.$\ Then, from (6.15), (6.11)-(i), (6.17) and Lemmas 2.2
(comparison theorem) and 2.3 we obtain, \be
\begin{array}{llll}
W_\delta(t,x)&\geq&  G^{t,x;u_1,
\beta^\varepsilon_1(u_1)}_{t,t+\delta} [W(t+\delta, X^{t,x;u_1,
\beta_1^\varepsilon(u_1)}_{t+\delta})]-\varepsilon\\
&\geq&G^{t,x;u_1,
\beta^\varepsilon_1(u_1)}_{t,t+\delta}[W(t+\delta, [X^{t,x;u_1,
\beta_1^\varepsilon(u_1)}_{t+\delta}])-C\varepsilon]-\varepsilon\\
&\geq&G^{t,x;u_1,
\beta^\varepsilon_1(u_1)}_{t,t+\delta}[W(t+\delta, [X^{t,x;u_1,
\beta_1^\varepsilon(u_1)}_{t+\delta}])]- C\varepsilon\\
&=&G^{t,x;u_1,\beta_1^\varepsilon(u_1)}_{t,t+\delta}[\sum\limits_{i\geq1}\textbf{1}_{\{X^{t,x;u_1,
\beta_1^\varepsilon(u_1)}_{t+\delta}\in O_i\}}W(t+\delta,y_i)]-
C\varepsilon,\ \ \mbox{P-a.s.}
\end{array}
\ee Furthermore, from (6.18), (6.11)-(ii), (6.16) and Lemmas 2.2 (comparison theorem) and 2.3, we have, \be
\begin{array}{lcl}
W_\delta(t,x)&\geq&
G^{t,x;u_1,\beta_1^\varepsilon(u_1)}_{t,t+\delta}[\sum\limits_{i\geq1}\textbf{1}_{\{X^{t,x;u_1,
\beta_1^\varepsilon(u_1)}_{t+\delta}\in O_i\}}J(t+\delta, y_i;
u_2, \beta^\varepsilon_{y_i}(u_2))-\varepsilon]-
C\varepsilon\\
&\geq&G^{t,x;u_1,\beta_1^\varepsilon(u_1)}_{t,t+\delta}[\sum\limits_{i\geq1}\textbf{1}_{\{X^{t,x;u_1,
\beta_1^\varepsilon(u_1)}_{t+\delta}\in O_i\}}J(t+\delta, y_i;
u_2, \beta^\varepsilon_{y_i}(u_2))]-C\varepsilon\\
&=&G^{t,x;u_1,\beta_1^\varepsilon(u_1)}_{t,t+\delta}[J(t+\delta,
[X^{t,x;u_1, \beta_1^\varepsilon(u_1)}_{t+\delta}]; u_2,
\beta^\varepsilon_{u_1}(u_2))]-C\varepsilon\\
&\geq &G^{t,x;u_1,
\beta^\varepsilon_1(u_1)}_{t,t+\delta}[J(t+\delta,X^{t,x;u_1,
\beta_1^\varepsilon(u_1)}_{t+\delta}; u_2,
\beta^\varepsilon_{u_1}(u_2))-C\varepsilon]- C\varepsilon\\
&\geq &G^{t,x;u_1,
\beta^\varepsilon_1(u_1)}_{t,t+\delta}[J(t+\delta,X^{t,x;u_1,
\beta_1^\varepsilon(u_1)}_{t+\delta}; u_2,
\beta^\varepsilon_{u_1}(u_2))]- C\varepsilon\\
&=& G^{t,x;u,\beta^\varepsilon(u)}_{t,t+\delta}[Y_{t+\delta}^{t,
x, u, \beta^{\varepsilon}(u)}]- C\varepsilon\\
&=& Y_{t}^{t, x; u, \beta^{\varepsilon}(u)}- C\varepsilon,\
\mbox{P-a.s., for any}\ u\in {\mathcal{U}}_{t, T}.
\end{array}
\ee Consequently, \be
\begin{array}{llll}
W_\delta(t,x)&\geq& \mbox{esssup}_{u \in {\mathcal{U}}_{t,
T}}J(t, x; u, \beta^{\varepsilon}(u))- C\varepsilon\\
&\geq&\mbox{essinf}_{\beta \in {\mathcal{B}}_{t,
T}}\mbox{esssup}_{u \in {\mathcal{U}}_{t, T}}J(t, x; u,
\beta(u))- C\varepsilon\\
&=&W(t,x)- C\varepsilon,\ \mbox{P-a.s.}
\end{array}
\ee Finally, letting $\varepsilon\downarrow0$\ we get
$W_\delta(t,x)\geq W(t,x).$\ The proof is complete.\endpf

\br\mbox{}{\rm{(i)}} From the inequalities (6.10) and (6.15) we see
that for all $(t, x)\in [0,T]\times {\mathbb{R}}^n,$\ $\delta>0$\
with $\delta\leq T-t$\ and $\varepsilon>0$,\ it
holds:\\
 a) For every $\beta \in {\cal{B}}_{t, t+\delta},$\
there exists some $u^{\varepsilon}(\cdot) \in {\cal{U}}_{t,
t+\delta}$\ such that
 \be W(t,x)(=W_\delta(t, x))\leq G^{t,x;
u^{\varepsilon},\beta(u^{\varepsilon})}_{t,t+\delta}
      [W(t+\delta, X^{t,x; u^{\varepsilon},\beta(u^{\varepsilon})}_{t+\delta})]+
      \varepsilon,\ \mbox{P-a.s.}
\ee
 b) There exists some $\beta^{\varepsilon} \in {\cal{B}}_{t,
t+\delta}$\ such that, for all $u\in {\cal{U}}_{t, t+\delta},$ \be
W(t,x)(=W_\delta(t, x))\geq G^{t,x;
u,\beta^{\varepsilon}(u)}_{t,t+\delta}
      [W(t+\delta, X^{t,x;u,\beta^{\varepsilon}(u)}_{t+\delta})]-
      \varepsilon,\ \mbox{P-a.s.}
\ee {\rm{(ii)}} Recall that the lower value function $W$\ is
deterministic. Thus, by choosing $\delta=T-t$\ and taking the
expectation on both sides of (6.21) and (6.22) we can show that $$
W(t,x)= \mbox{inf}_{\beta \in {\cal{B}}_{t,T}}\mbox{sup}_{u \in
{\mathcal{U}}_{t,T}}E[J(t,x; u,\beta(u))]. $$ In analogy we also
have
$$ U(t,x)= \mbox{sup}_{\alpha
\in {\cal{A}}_{t,T}}\mbox{inf}_{v \in {\mathcal{V}}_{t,T}}E[J(t,x;
\alpha(v), v)]. $$ \er

\noindent{\bf Acknowledgment}\ Juan Li thanks the Department of
Mathematics of the University of West Brittany, and especially
Rainer Buckdahn, for their hospitality during her stay in France.

\end{document}